\documentclass[11pt]{article}
\usepackage{amsfonts}
\newcommand{\be}{\begin{equation}}
\newcommand{\ee}{\end{equation}}
\newcommand{\bea}{\begin{eqnarray}}
\newcommand{\eea}{\end{eqnarray}}
\newcommand{\beaa}{\begin{eqnarray*}}
\newcommand{\eeaa}{\end{eqnarray*}}
\newcommand{\eq}[1]{(\ref{#1})}

\newcommand{\eps}{\varepsilon}

\newcommand{\metric}{\mathbf{g}}
\newcommand{\zero}[1]{\stackrel{0}{#1}}
\newcommand{\one}[1]{\stackrel{1}{#1}}
\newcommand{\two}[1]{\stackrel{2}{#1}}
\newcommand{\mone}[1]{\stackrel{-1}{#1}}
\newcommand{\mtwo}[1]{\stackrel{-2}{#1}}
\newcommand{\mthree}[1]{\stackrel{-3}{#1}}
\newcommand{\mfour}[1]{\stackrel{-4}{#1}}

\def\Rl{\mathbb{R}}
\def\Deta{D_\eta}
\def\Tbar{{T}}
\def\Ybar{{Y}}
\def\Sone{S^{(1)}}
\def\Lmtwo{L^{(-2)}}

\newtheorem{proposition}{Proposition}
\newtheorem{definition}{Definition}

\begin{document}

\renewcommand{\theequation}{\thesection.\arabic{equation}}

\title{Diffractive
Nonlinear Geometrical Optics \\ for Variational
Wave Equations \\
and the Einstein Equations}
\author{G.~Al\`{\i} \\ Istituto per le Applicazioni del Calcolo \\
Consiglio Nazionale delle Ricerche, Napoli \\
and INFN-Gruppo c. Cosenza
\and J.~K.~Hunter
\\ Department of Mathematics \\
University of California at Davis}

\date{October 31, 2005}

\maketitle

\begin{abstract}
We derive an asymptotic solution of the
vacuum Einstein equations that describes
the propagation and diffraction of a localized, large-amplitude,
rapidly-varying gravitational wave.
We compare and contrast the resulting theory of strongly
nonlinear geometrical optics for the Einstein equations
with nonlinear geometrical optics theories for
variational wave equations.
\end{abstract}

\section{Introduction}
\setcounter{equation}{0}

Geometrical optics\footnote{Here, we use the term
`geometrical optics' to refer to any asymptotic theory for the propagation
of short-wavelength, high-frequency waves, irrespective of its area of physical application.}
and its generalizations, such as the geometrical theory of diffraction,
are a powerful approach to the study of wave propagation,
for both linear and nonlinear waves.
In this paper, we develop a theory of strongly nonlinear
geometrical optics for gravitational wave
solutions of the vacuum Einstein equations.
Specifically, we derive asymptotic equations that describe
the diffraction of large-amplitude, rapidly-varying gravitational waves.

These equations are a generalization of the straightforward
non-diffractive, nonlinear geometrical optics equations for large-amplitude
gravitational waves derived in \cite{AH}. This
strongly nonlinear theory differs fundamentally from the weakly
nonlinear theory for small-amplitude gravitational waves
obtained by Choquet-Bruhat \cite{CB} and Isaacson \cite{Is},
because it captures the direct nonlinear self-interaction of the waves.

The Einstein equations may be derived from a variational principle,
and, when written with respect to a suitable gauge,
they form a system of wave equations for the gravitational field.
In order to explain the structure
of nonlinear geometrical optics theories for gravitational waves,
and to motivate the form of our asymptotic expansion,
it is useful to consider first such theories for
a general class of variational wave equations.

We describe straightforward and diffractive
geometrical optics theories for variational wave equations in Section~\ref{ngo}.
In the weakly nonlinear theory for waves with periodic waveforms,
the amplitude-waveform function $a(\theta,v)$ of the wave
depends on a ray variable $v$ and periodically
on a `fast' phase variable $\theta$. The wave-amplitude $a$
satisfies the Hunter-Saxton equation \eq{HSeq},
\be
\left\{a_v + \left(\frac{1}{2}\Lambda a^2\right)_\theta + N a\right\}_\theta
= \frac{1}{2}\Lambda\left\{ a_\theta^2 - \langle a_\theta^2\rangle\right\},
\label{HSeq1}
\ee
where the angular brackets denote an average with respect to $\theta$.
The coefficient $\Lambda$ of the nonlinear terms in \eq{HSeq1} may be interpreted
as a derivative of the wave speed with respect to the wave amplitude.
In addition, the wave-amplitude is coupled with a slowly-varying mean field.

An application of this expansion to the Einstein equations leads to the
theory of Choquet-Bruhat \cite{CB} and Isaacson \cite{Is}.
In that case,
the nonlinear coefficients corresponding to $\Lambda$ in the equations
for the amplitudes of the gravitational waves are identically zero.
This theory therefore describes a nonlinear interaction
between a high-frequency, oscillatory gravitational wave and
a slowly-varying mean gravitational field, but does not describe the direct
nonlinear self-interaction of the gravitational wave itself.

In Section~\ref{lindegwave}, we make a distinction between `genuinely
nonlinear' wave-fields in variational wave equations, for which $\Lambda$ is never zero,
and `linearly degenerate' wave-fields
for which $\Lambda$ identically zero (see Definition~\ref{def:lindeg}).
All wave-fields in the Einstein equations are linearly degenerate. As
observed in \cite{CB2, CB3,GHZ,Hu,lindblad}, for example, this fact
reflects a fundamental degeneracy in the nonlinearity of the Einstein equations
in comparison with general wave equations. 

A property of the Einstein equations related to their linear degeneracy
is that they possess an exact solution for
non-distorting, large-amplitude, plane waves,
the Brinkman solution \cite{Br,LL}. One can use this
solution to derive a strongly nonlinear geometrical optics theory for
large-amplitude gravitational waves. We outline the resulting non-diffractive
theory \cite{AH} in Section~\ref{collidingpw}.

We derive our diffractive, strongly nonlinear geometrical optics solution
of the vacuum Einstein equations in Section~\ref{plane}. This
solution describes the propagation and diffraction
of a thin, large-amplitude gravitational
wave, such as a pulse or `sandwich' wave.
The simplest, and basic, case is that of a plane-polarized
gravitational wave diffracting in a single direction. We
summarize the resulting asymptotic equations here.

We suppose that the polarization
of the wave is aligned with the diffraction direction.
Then, with respect to a suitable coordinate system $(u,v,y,z)$,
the metric $\metric$ of the wave adopts the form
\be
\metric = -2e^{-M}\bigl(du - \eps \Ybar dy \bigr)dv
+e^{-U}\bigl(e^{V}dy^2+e^{-V}dz^2 \bigr) +O(\eps^2).
\label{metric}
\ee
Here, $\eps$ is a small parameter. The leading-order metric component
functions $(U, V, M)$, and the first-order function $\Ybar$
depend upon a `fast' phase variable $\theta$,
an `intermediate' transverse variable $\eta$,
and the `slow' ray variable $v$, where
\[
\theta = \frac{u}{\eps^2},\qquad \eta = \frac{y}{\eps}.
\]
The phase $u$ is light-like and the
transverse coordinate $y$ is space-like.

This metric describes a gravitational wave whose wavefronts are close to
the null-hypersurface $u=0$. The wave is plane-polarized in the
$(y,z)$-directions and diffracts in the $y$-direction.

To write equations for the metric component functions in a concise form, we
define a derivative $\Deta$ and functions $\phi$, $\psi$ by
\bea
\Deta &=& e^{U}\left(\partial_\eta + \Ybar\partial_\theta\right),
\label{defdeta}\\
\phi &=& \Deta M - e^{U}\Ybar_\theta,
\label{defphi}\\
\psi &=& \Deta(U+V).
\label{defpsi}
\eea
We note that $\Deta$ has the geometrical significance that
\[
\Deta = \frac{1}{\eps} e^V (dy)^\sharp,
\]
where $\sharp$ denotes the `raising' operator from one-forms to vector fields.

Then $(U, V, M, \Ybar)$ satisfy the following system of PDEs:
\bea
&&U_{\theta\theta} - \frac{1}{2}\left(U_\theta^2 + V_\theta^2\right) + U_\theta M_\theta = 0,
\label{eq1}\\
&&(\phi+\psi)_\theta = \psi (U+V)_\theta,
\label{eq2}\\
&&
U_{\theta v}- U_{\theta}U_{v}
=\frac12 e^{-(U+V+M)}\left\{\Deta\phi+\Deta\psi
-\frac12 \phi^2 - \phi\psi - \psi^2\right\},
\label{eq3}
\\
&&
V_{\theta v}
-\frac{1}{2}\left( U_{\theta}V_{v}+U_{v}V_{\theta}\right)
=\frac12 e^{-(U+V+M)}\left\{-\Deta\phi+\frac12 \phi^2\right\},
\label{eq4}
\\
&&
M_{\theta v}
+\frac{1}{2}\left( U_{\theta}U_{v}-V_{\theta}V_{v}\right)
=\frac12 e^{-(U+V+M)}\left\{-\Deta\psi
-\frac12 \phi^2 + \psi^2\right\}.
\qquad
\label{eq5}
\eea
This system is the main result of our analysis. It is an asymptotic reduction of the
full vacuum Einstein equations to a $(1+2)$-dimensional system of PDEs. The system
provides a model nonlinear wave equation for general relativity, and should
be useful, for example, in studying the focusing of nonplanar gravitational
waves and the effect of diffraction on the formation of singularities. We plan
to study these topics in future work.

We also derive
asymptotic equations for the diffraction of a gravitational wave of general
polarization in two transverse directions. These equations are given by
\eq{thetaconstraint}--\eq{apc1}, \eq{Req1}, and \eq{R0a-0}--\eq{Rab-1},
but they
are much more complicated than the ones for plane-polarized waves written out above.

We study the structure of \eq{eq1}--\eq{eq5} in Section~\ref{properties}.
The constraint equation \eq{eq1} is a nonlinear
ODE in $\theta$ in which $\eta$, $v$ appear as parameters.
As shown in Proposition~\ref{prop:con}, this constraint
is preserved by the remaining equations. Equation \eq{eq2} is a linear,
nonhomogeneous ODE in $\theta$ for
$\Ybar$, whose coefficients depend on $\theta$- and
$\eta$-derivatives of $(U,V,M)$. It may therefore be regarded as determining
$Y$ in terms of $(U,V,M)$.
Equations \eq{eq3}--\eq{eq5} are a system of
evolution equations for $(U,V,M)$ which are coupled with $\Ybar$.
The main part of the system
of evolution equations consists of a (1+2)-dimensional wave equation
in $(\theta,\eta,v)$ for $(U+V)$, in which $(\theta,v)$ are characteristic
coordinates.

When $\Ybar=0$ and all functions are independent of $\eta$, equations \eq{eq1}--\eq{eq5} reduce to the
colliding plane wave equations, without the $v$-constraint equation (see Section~\ref{collidingpw}).
When all functions depend on $(\xi,\eta)$ with $\xi = \theta - \lambda v$
for some constant $\lambda$, we get a system of PDEs in two variables that is studied further
in \cite{abh}. This system describes space-times that are stationary with respect to an
observer moving close to the speed of light.

\section{Geometrical optics}
\label{ngo}
\setcounter{equation}{0}

The Einstein equations do not form a hyperbolic system of PDEs
because of their gauge-covariance, but their properties
are closely related to those of hyperbolic PDEs.
In order to develop and interpret geometrical optics solutions for the Einstein
equations, it is useful to begin by studying geometrical optics solutions
for hyperbolic systems of variational wave equations.

We remark that there is a close correspondence between geometrical optics
theories for hyperbolic systems of conservation
laws \cite{Hrev} and variational wave equations. For example, the
inviscid Burgers equation \cite{CB1,HK} is the analog of the Hunter-Saxton equation
\eq{hsloc}, and the unsteady transonic small disturbance equation \cite{Hu} is the
analog two-dimensional Hunter-Saxton equation \eq{diffHS}. One can also derive
large-amplitude geometrical optics theories for linearly degenerate waves
in hyperbolic conservation laws (see \cite{serre1,serre2}, for example)
that are analogous to the large-amplitude theories described here for the
Einstein equations. 

There are other nonlinear geometrical optics theories for dispersive waves,
most notably Whitham's `averaged Lagrangian method' \cite{taub,Wh} for large-amplitude dispersive
waves. This theory has a different character
from the ones for nondispersive hyperbolic waves.
Nonlinear dispersive waves have specific waveforms ---
given by traveling wave solutions --- in which the effects of dispersion and nonlinearity balance,
whereas nondispersive hyperbolic wave equations,
and the Einstein equations, have traveling wave solutions
with arbitrary waveforms, which may distort as the
waves propagate.

In straightforward theories of geometrical optics, waves
are locally approximated by plane waves and propagate along rays. In linear geometrical optics, the wave
amplitude satisfies an ODE (the transport equation) along a ray; in nonlinear
geometrical optics, the amplitude-waveform function typically satisfies a nonlinear PDE in one space
dimension (a generalization of the transport equation) along a ray.
This difference may be understood as follows:
linear hyperbolic waves propagate without distortion,
so one requires only an ODE along each ray to determine
the change in the wave amplitude; by contrast, wave-steepening and other effects
typically distort the waveform of nonlinear hyperbolic waves,
so one requires a PDE along each ray to determine the change in the wave amplitude
and the waveform.

Straightforward geometrical optics, and a local plane-wave approximation,
break down when the effects of wave-diffraction become important; for example,
this occurs when a high-frequency wave focuses at a caustic, or when a wave
beam of large (relative to its wavelength), but finite, transverse extent spreads
out. The effects of diffraction on a high-frequency wave may be described by
the inclusion of additional length-scales in the straightforward geometrical
optics asymptotic solution.

Perhaps the most basic asymptotic solution that incorporates the effect
of wave diffraction leads to the `parabolic approximation',
described below for wave equations.
In Section~\ref{plane}, we will
derive analogous asymptotic solutions of the Einstein equations
for the diffraction of large-amplitude gravitational waves.

\subsection{The wave equation}
\label{sec:wave}

We begin by recalling geometrical optics theories for the linear wave equation
\be
g_{tt} = \nabla \cdot\left(c_0^2 \nabla g\right).
\label{linwave}
\ee
Here, $g(t,x)$ is a scalar function, the wave speed $c_0(t,x)$ is a given
smooth function, and $x\in \Rl^d$. Although well-known, these theories
provide a useful background to our analysis of nonlinear variational wave equations
and the Einstein equations.

We look for a short-wavelength asymptotic solution $g = g^\eps$
of \eq{linwave}, depending on a small parameter $\eps$, of the form
\[
g^\eps(t,x) \sim a\left(\frac{u(t,x)}{\eps^2},t,x\right)
\qquad \mbox{as $\eps\to 0$}.
\]
The solution depends on a `fast' phase variable\footnote{Here, we use
$u/\eps^2$ as a phase variable, rather than $u/\eps$,
for consistency with the diffractive expansion below.}
\[
\theta=\frac{u}{\eps^2}
\]
and the `slow' space-time variables $(t,x)$.

One finds that the scalar-valued phase function $u(t,x)$ satisfies
the eikonal equation
\be
u_t^2 = c_0^2 |\nabla u|^2. \label{lineik}
\ee
The amplitude-waveform function $a(\theta,t,x)$ satisfies the equation
\be
a_{\theta v} + N a_\theta = 0, \label{wavetran} \ee
where
\be
\partial_v = u_t\partial_t - c_0^2{\nabla u}\cdot\nabla
\label{linray}
\ee
is a derivative along the rays associated with the phase $u$,
and $N(t,x)$ is given by
\be
N = \frac{1}{2}\left\{u_{tt} - \nabla\cdot\left(c_0^2\nabla u\right)\right\}.
\label{defN}
\ee

Equation \eq{wavetran} has solutions of the form
\[
a(\theta,t,x) =  A(t,x) F\left(\theta\right),
\]
where $A(t,x)$ is a wave-amplitude, which satisfies an ODE
along a ray (the transport equation of linear geometrical optics)
\be
A_v + N A = 0,
\label{ODEte}
\ee
and $F(\theta)$ is an arbitrary function that
describes the waveform of the wave.
For example, if $F(\theta) = e^{i\theta}$, then the
solution describes an oscillatory harmonic wave; if
\[
F(\theta) = \left\{ \begin{array}{rl}
\theta^n & \mbox{$\theta > 0$},
\\
0 & \mbox{$\theta \le 0$},
\end{array}\right.
\]
then the solution describes a wavefront across which the normal derivative of $g$ of order
$n$ jumps; and if $F(\theta)$ has compact support, then the solution describes a localized pulse.

The term $N A$ in the transport equation \eq{ODEte} describes the effect of the ray geometry on the wave
amplitude. The coefficient $N$ becomes infinite on caustics where the rays focus.
The straightforward geometrical optics solution breaks down when this happens, and diffractive effects
must then be taken into account.

There are many ways in which diffraction modifies straightforward
geometrical optics. Here, we consider one of the simplest diffractive
expansions, given by
\[
g^\eps(t,x) \sim a\left(\frac{u(t,x)}{\eps^2},\frac{y(t,x)}{\eps},t,x\right)
\qquad \mbox{as $\eps\to 0$}.
\]
This asymptotic solution depends upon an additional `intermediate' variable
\[
\eta = \frac{y}{\eps},
\]
where $y(t,x)$ is a scalar-valued transverse phase.

One finds that $u$ satisfies the eikonal equation, as before, and
\[
y_v = 0,
\]
meaning that
$y$ is constant along the rays associated with $u$. Moreover, the amplitude-waveform
function $a(\theta,\eta,t,x)$ satisfies the equation
\be
a_{\theta v}+ N a_\theta  +  \frac{1}{2} D a_{\eta\eta} = 0,
\label{genpareq}
\ee
where the coefficient $D(t,x)$ of the diffractive term is given by
\be
D = y_t^2 - c_0^2|\nabla y|^2.
\label{defDD}
\ee
For harmonic solutions, we have
\[
a(\theta,\eta,t,x) = A(\eta,t,x) e^{i\theta},
\]
and equation \eq{genpareq} reduces to a Schr\"odinger equation
\[
i \left\{A_v + N A\right\}+ \frac{1}{2} D A_{\eta\eta}  = 0.
\]
This `parabolic approximation' and its generalizations 
are widely used in the study of wave propagation \cite{parabolic}.

In the simplest case of the diffraction of plane wave
solutions of the two-dimensional wave equation with
wave speed $c_0=1$,
\[
g_{tt} = g_{xx} + g_{yy},
\]
we may choose
\[
u = \frac{t-x}{\sqrt{2}},\qquad y =y, \qquad v = \frac{t+x}{\sqrt{2}}.
\]
Equation \eq{genpareq} is then
\be
a_{\theta v} = \frac{1}{2}a_{\eta\eta}.
\label{pareq}
\ee
This equation describes waves that
propagate in directions close to the $x$-direction, and is a wave
equation in which $\theta$ and $v$ are characteristic coordinates.
We will see that equations with a similar structure to \eq{pareq} in their
highest-order derivatives arise from the Einstein equations.

\subsection{A variational wave equation}

Next, we consider the following nonlinear, scalar wave equation \cite{HS}
\be
g_{tt} - \nabla\cdot \left(c^2(g)\nabla g\right) +
c(g) c^\prime(g) |\nabla g|^2 = 0.
\label{nonwave}
\ee
We assume that the wave speed $c:\Rl\to\Rl^+$ is a smooth,
non-vanishing function, and a prime denotes the derivative with respect to $g$.
This equation is derived from the variational principle
\[
\delta\int \left\{\frac{1}{2} g_t^2 - \frac{1}{2}c^2(g) |\nabla g|^2\right\}\, dt dx = 0.
\]

The structure of the nonlinear terms in \eq{nonwave} resembles that of the
Einstein equations, although, as we shall see, the effects of nonlinearity
are qualitatively different because of a `linear degeneracy' in the Einstein equations.

We look for an asymptotic solution of \eq{nonwave} of the form \cite{GHZ,Hrev,HS}
\[
g^\eps(t,x) \sim g_0(t,x) +
\eps^2 a\left(\frac{u(t,x)}{\eps^2},t,x\right)
\qquad\mbox{as $\eps \to 0$.}
\]
This solution represents a small-amplitude, high-frequency perturbation
of a slowly-varying field $g_0$. The amplitude and the wavelength
of the perturbation
are chosen to be of the same order of magnitude because this leads to
a balance between the effects of weak nonlinearity and the ray geometry.

First, we suppose that the amplitude-waveform function $a(\theta,t,x)$
is a periodic function of the phase variable $\theta$. We assume, without
loss of generality, that its mean
with respect to $\theta$ is zero.

We find  that the phase $u(t,x)$ satisfies the linearized eikonal equation
\eq{lineik}, with $c_0 = c\left(g_0\right)$.
The mean-field $g_0$ satisfies the nonlinear wave equation
\be
g_{0tt} - \nabla\cdot\left(c_0^2 \nabla g_0\right) + c_0 c_0^\prime |\nabla g_0|^2
+ |\nabla u|^2 c_0 c_0^\prime \langle a^2_\theta\rangle = 0,
\label{wavemean}
\ee
where the angular brackets denote an average with respect to $\theta$
over a period, and
$
c_0^\prime = c^\prime\left(g_0\right)
$.
Equation \eq{wavemean} has the same form as the original wave equation, with
an additional source term proportional to the mean energy of the wave-field.

The amplitude-waveform function $a(\theta,t,x)$ satisfies the periodic Hunter-Saxton equation
\eq{HSeq1},
\be
\left\{a_v + \left(\frac{1}{2}\Lambda a^2\right)_\theta + N a\right\}_\theta
= \frac{1}{2}\Lambda\left\{ a_\theta^2 - \langle a_\theta^2\rangle\right\},
\label{HSeq}
\ee
where $\partial_v$ is the ray derivative defined in \eq{linray}, $N$ is given by
\eq{defN}, and
\be
\Lambda = - |\nabla u|^2 c_0 c_0^\prime.
\label{defLambdawave}
\ee
The coefficient $\Lambda$ is proportional to the derivative of the wave-speed
with respect to the wave amplitude, and it provides crucial information about the effect
of nonlinearity on the waves; for a given $g_0$,
we obtain a nonlinear PDE for the amplitude-waveform
function only if $\Lambda\ne 0$.

The expansion summarized above for waves with periodic waveforms is uniformly valid for $(t,x)=O(1)$
and $u = O(\eps^{-1})$ as $\eps \to 0$.
If we consider localized waves, then we may obtain an asymptotic expansion
that is valid near\footnote{Specifically, the expansion is valid when $\theta = O(1)$,
or $u=O(\eps^2)$, and $t=O(1)$.}  the wave-front $u=0$, in which we neglect the mean-field effects.
Thus, the background field $g_0(x)$ is a solution
of the original wave equation, and $a(\theta,t,x)$ satisfies
\be
\left\{a_v + \left(\frac{1}{2}\Lambda a^2\right)_\theta + N a\right\}_\theta
= \frac{1}{2}\Lambda a_\theta^2.
\label{hsloc}
\ee

It is straightforward to include diffractive effects in the expansion
for localized waves. The asymptotic solution has the form
\[
g^\eps(t,x) \sim g_0(x) +
\eps^2 a\left(\frac{u(t,x)}{\eps^2},\frac{y(t,x)}{\eps},t,x\right)
\qquad\mbox{as $\eps \to 0$.}
\]
Here, $g_0$ is a solution of the original wave equation,
the phase $u$ satisfies the eikonal equation at $g_0$,
the transverse function $y$ is constant
along the rays associated with $u$, and
$a(\theta,\eta, t,x)$ satisfies a $(1+2)$-dimensional generalization
of the Hunter-Saxton equation \eq{hsloc},
\be
\left\{a_v + \left(\frac{1}{2}\Lambda a^2\right)_\theta + N a\right\}_\theta
+\frac{1}{2} D a_{\eta\eta} = \frac{1}{2}\Lambda a_\theta^2,
\label{diffHS}
\ee
where $D$ is given by \eq{defDD}.

\subsection{Systems of variational wave equations}
\label{sysvar}

In this section, we consider a class of hyperbolic systems
of nonlinear wave equations that are derived from variational principles of the form \cite{HS}
\be
\delta\int A^{\alpha\beta}_{pq}(g){\frac{\partial g^p}{\partial x^\alpha}}
\frac{\partial g^q}{\partial x^\beta}\,dx  = 0.
\label{(B.1)}
\ee
Here,
$x = (x^0,\dots, x^d)\in\Rl^{d+1}$ are the space-time variables,
\[
g = \left(g^1,\dots , g^m\right) : {\Rl}^{d+1}\to{\Rl}^m
\]
are the dependent variables,
$
A^{\alpha\beta}_{pq} : \Rl^m \to \Rl
$
are smooth coefficient functions, and we use the summation convention.
We assume that
$
A^{\alpha\beta}_{pq} = A^{\beta\alpha}_{pq} = A^{\alpha\beta}_{qp}
$.

The Euler-Lagrange equations associated with \eq{(B.1)} are
\be
G_p\left[g\right] = 0,
\label{(B.2)}
\ee
where
\be
G_p[g] = \frac{\partial}{\partial x^\alpha}\left\{A^{\alpha\beta}_{pq}\left(g\right)
\frac{\partial g^q}{\partial x^\beta}\right\}
-\frac12\frac{\partial A^{\alpha\beta}_{qr}}{\partial g^p}\left(g\right)
\frac{\partial g^q}{\partial x^\alpha}\frac{\partial g^r}{\partial x^\beta}.
\label{defG}
\ee
We assume that \eq{(B.2)} forms a hyperbolic system of PDEs.
The scalar wave equation considered in the previous section is the
simplest representative of this class of equations.

The weakly nonlinear geometrical optics solution of \eq{(B.2)}
has the form \cite{GHZ,Hrev,HS}
\be
g^\eps(x) \sim   g_0(x) + \eps^2  a\left(\frac{u(x)}{\eps^2},x\right)R(x) \qquad \mbox{as $\eps\to 0$},
\label{(B.4)}
\ee
where $g_0 : \Rl^{d+1} \to \Rl^m$ is a slowly-varying function,
$u : \Rl^{d+1}\to \Rl$ is a phase function, $a: \Rl\times\Rl^{d+1}\to \Rl$
is an amplitude-waveform function, and $R: \Rl^{d+1}\to \Rl^m$ is a suitable vector field.
We find that the phase $u$ satisfies the eikonal equation
\be
\det C = 0,
\label{(B.3)}
\ee
where the $m\times m$ matrix $C(x)$ has components
\be
C_{pq} = u_{x^\alpha}u_{x^\beta} A^{\alpha\beta}_{pq}\left(g_0\right).
\label{defmatC}
\ee
The vector $R = (R^1,\dots,R^m)^T$ in \eq{(B.4)} is a right null-vector
of $C$, so that $C_{pq}R^q = 0$. Here, and in \eq{(B.4)},
we assume that we are dealing with a simple characteristic, meaning
that the null-space of $C$ is one-dimensional.

We write the scalar-valued amplitude-waveform function $a(\theta,x)$
as a function of the `fast' phase variable $\theta = u/\eps^2$,
and the slow variables $x$. If $a(\theta,x)$ is a periodic
function of $\theta$ with zero mean, then the mean-field
$g_0(x)$ satisfies the equation
\be
G_p\left[g_0\right] = \frac{1}{2} H_p\left\langle a_\theta^2\right\rangle,
\label{mfeq}
\ee
where the angular brackets denote an average with respect to $\theta$,
and
\[
H_p  = u_{x^\alpha}u_{x^\beta}\frac{\partial A^{\alpha\beta}_{qr}}{\partial g^p}\left(g_0\right)R^q R^r.
\]
Equation \eq{mfeq} has the same form as the original equation for $g$
with an additional source term proportional to $\langle a_\theta^2\rangle$.

The equation for $a$ is \eq{HSeq} with
\bea
\partial_v &=& 2u_{x^\beta} A^{\alpha\beta}_{pq}\left(g_0\right)R^pR^q \partial_{x^\alpha},
\nonumber\\
\Lambda &=& u_{x^\alpha}u_{x^\beta}\frac{\partial A^{\alpha\beta}_{qr}}{\partial g^p}\left(g_0\right)R^p R^q R^r,
\label{defvarLambda}\\
N &=& \frac{\partial}{\partial x^\alpha}\left\{u_{x^\beta} A^{\alpha\beta}_{pq}\left(g_0\right)R^p R^q\right\}
- u_{x^\alpha} \frac{\partial A^{\alpha\beta}_{qr}}{\partial g^p}
\left(g_0\right)\frac{\partial g_0^r}{\partial x^\beta}R^p R^q.
\nonumber
\eea
For localized solutions, the mean-field interactions may be neglected, and we can obtain diffractive
versions of this expansion as before.

These equations generalize easily to the case when equation \eq{(B.2)}
has a multiple characteristic of constant multiplicity $n \ge 2$, say. In that case,
one obtains a mean-field equation whose source term is a sum
of averages of products of $\theta$-derivatives of the wave amplitudes,
and an $n\times n$ coupled system of Hunter-Saxton equations for the
wave amplitudes. In particular, if
the vectors $\{R_1,\dots,R_n\}$ form a
basis of the null-space of the matrix $C$
defined in \eq{defmatC}, then the coefficients
$\Lambda_{ijk}$ of the nonlinear terms in the
system of equations for the wave amplitudes are given by
\be
\Lambda_{ijk} =  u_{x^\alpha}u_{x^\beta}\frac{\partial A^{\alpha\beta}_{qr}}{\partial g^p}
\left(g_0\right)R_i^p R_j^q R_k^r,
\qquad 1\le i,j,k \le n.
\label{defvarLambdamult}
\ee

\subsection{Linearly degenerate wave equations}
\label{lindegwave}

Motivated by the corresponding definition for hyperbolic conservation
laws introduced by Lax \cite{lax}, we make the following definition.

\begin{definition}
\label{def:lindeg}
A simple characteristic of a hyperbolic system of wave equations
\eq{(B.2)}--\eq{defG} is \emph{genuinely
nonlinear} (respectively, \emph{linearly degenerate})
if, for every $g_0\in \Rl^m$ and every non-zero $du\in \Rl^{d+1}$
that satisfies \eq{(B.3)}, the quantity
$\Lambda$ defined in \eq{defvarLambda} is non-zero
(respectively, zero). A multiple characteristic
of constant multiplicity $n$ is \emph{linearly degenerate}
if all coefficients $\Lambda_{ijk}$ defined in \eq{defvarLambdamult} are zero.
We say that the system is \emph{genuinely nonlinear} (respectively, \emph{linearly degenerate})
if all of its characteristics are genuinely nonlinear (respectively, linearly degenerate).
\end{definition}

Thus, for linearly degenerate wave equations, the amplitude-waveform function in the weakly
nonlinear theory corresponding to the ansatz \eq{(B.4)} always satisfies
a linear PDE. Interactions with a mean-field may still occur, however.

For example, from \eq{defLambdawave}, the nonlinear wave equation \eq{nonwave}
is genuinely nonlinear if $c^\prime(g) \ne 0$ for all $g\in \Rl$, and
linearly degenerate if $c$ is constant, when it reduces to the linear wave equation.

Even if a wave-field is not linearly degenerate, a loss of
genuine nonlinearity may occur at a particular point $g_0\in \Rl^m$ and direction
$du \in \Rl^{d+1}$. A loss of genuine nonlinearity
is related to the null condition for nonlinear wave equations
introduced by Klainerman  \cite{klainerman}. The Einstein equations do not satisfy the
Klainerman null condition completely because of the non-zero coupling between gravitational waves and a
mean gravitational field \cite{CB2}. Nevertheless, it is the vanishing of the nonlinear $\Lambda$-coefficients
in the equations for the amplitude-waveform functions at the Minkowski metric, for example,
that permits the existence of global smooth, small-amplitude perturbations
\cite{CK,nicolo,lindblad1}. A similar result would not be true for $(1+3)$-dimensional
variational wave equations in which $\Lambda\ne 0$.

We emphasize that this definition for
wave equations is analogous to but different from the definition for hyperbolic systems
of conservation laws. If a genuinely nonlinear wave equation
were rewritten as a first-order hyperbolic system, it
would be classified as a linearly degenerate first-order hyperbolic
system, since its wave speeds are independent of the derivative $\nabla g$.

Hyperbolic systems of conservation laws with linearly degenerate wave-fields
possess non-distorting plane wave solutions, and one can derive
a large-amplitude geometrical optics theory for them \cite{serre1,serre2}.
Linearly degenerate variational systems
of wave equations possess non-distorting plane wave solutions,
and a large-amplitude geometrical optics theory, only if
they satisfy certain additional degeneracy conditions.
We will not analyze these issues
here, however, and instead discuss the analogous
theory for the Einstein equations.

\subsection{The Einstein equations}
\label{collidingpw}

The Einstein equations do not fall exactly into the class of variational
wave equations considered in Section~\ref{sysvar} because of their gauge-covariance.
They can be derived from a variational principle of the form \eq{(B.1)}, obtained
after an integration by parts in the Einstein-Hilbert action, but the resulting
Euler-Lagrange equations are not hyperbolic. Nevertheless, as is well-known,
they become hyperbolic when written with respect to a suitable gauge, such as the harmonic gauge.
The equation for the metric $\metric$, with components $g_{\mu\nu}$ and determinant $g$, then adopts a
similar form to \eq{(B.2)}--\eq{defG}, namely
\[
\frac{\partial}{\partial x^\alpha}\left\{g^{\alpha\beta}\sqrt{-g}\,\frac{\partial g_{\mu\nu}}{\partial x^\beta} \right\}
- H_{\mu\nu}(\metric)\left(\partial \metric,\partial \metric\right) = 0,
\]
where $H_{\mu\nu}$ is a quadratic form in the metric derivatives
$\partial\metric$ with coefficients depending on $\metric$.
(See \cite{lindblad1}, for example, for an explicit expression.)

The weakly nonlinear expansion described in Section~\ref{sysvar}
for variational systems of wave equations corresponds
to the expansion of Choquet-Bruhat \cite{CB} and Isaacson \cite{Is} for the Einstein
equations. The Einstein equations are not strictly hyperbolic, and
all of the nonlinear coefficients in the asymptotic equations for the
amplitude-waveform functions are zero.
Thus, in this generalized sense, the Einstein equations are a linearly degenerate
system of variational wave equations.

The Einstein equations possess a non-distorting, plane wave
solution for large-amplitude gravitational waves,
which forms the basis of
a strongly nonlinear geometrical optics theory for
large-amplitude gravitational waves.
It does not, however, appear possible to obtain a self-consistent
theory for oscillatory large-amplitude waves,
since the mean energy-momentum associated with
an extended wave-packet would generate a
very strong background curvature of space-time.
This restriction is related to the fact that mean-field interactions
with oscillatory waves already occur at leading order in the small-amplitude theory.
We therefore consider localized waves.

An asymptotic theory for the propagation of localized, large-amplitude,
rapidly varying gravitational waves into slowly varying space-times
was developed in \cite{AH}. In this theory, the metric of a plane-polarized
wave may be written, with respect to a suitable coordinate
system $(u,v,y,z)$, as
\be
\metric=-2e^{-M}du\,dv+e^{-U}(e^{V}dy^2+e^{-V}dz^2) + O(\eps^2).\label{Br}
\ee
Here, the metric component functions $(U,V,M)$, depend
on $(\theta,v)$, where the `fast' phase variable $\theta$ is given by
\[
\theta = \frac{u}{\eps^2}.
\]
The functions $(U,V,M)$ satisfy the following PDEs:
\bea
&&U_{\theta \theta} - \frac{1}{2}\left(U_{\theta}^2+V_{\theta}^2\right)+U_{\theta}M_{\theta} = 0,
\label{pc1}
\\
&&U_{\theta v} - U_{\theta}U_{v} = 0,
\label{pc2}
\\
&&V_{\theta v} - \frac{1}{2}\left(U_{\theta}V_{v}+U_{v}V_{\theta}\right) = 0,
\label{pc3}
\\
&&M_{\theta v} + \frac{1}{2}\left(U_{\theta}U_{v}-V_{\theta}V_{v}\right) = 0.
\label{pc4}
\eea
Equations (\ref{pc2})--(\ref{pc4}) are wave equations
for $(U, V, M)$ in characteristic coordinates $(\theta,v)$,
and \eq{pc1} is a constraint which is preserved by (\ref{pc2})--(\ref{pc4}).

These equations correspond to a well-known exact solution of the Einstein equations,
the colliding plane wave solution \cite{Gr,KP,Sz1,Sz2}, without the
usual constraint equation in $v$,
\[
U_{v v} - \frac{1}{2}\left(U_{v}^2+V_{v}^2\right)+U_{v}M_{v} = 0.
\]
This constraint need not be satisfied if
the slowly-varying space-time into which the wave propagates
is not that of a counter-propagating gravitational wave. If it is not satisfied, then
the resulting metric is an asymptotic solution of the Einstein equations but not
an exact solution.

The asymptotic equations we derive in this paper
are a diffractive generalization
of (\ref{pc1})--(\ref{pc4}).

\section{The asymptotic expansion}
\label{plane}
\setcounter{equation}{0}

In this section, we outline our asymptotic expansion
of the Einstein equations, and specialize it to the case
of plane-polarized gravitational waves that diffract
in a single direction. Our goal is to explain
the structure of the expansion and the resulting perturbation
equations. The detailed algebra is summarized in the appendices.

\subsection{The general expansion}

Let $\metric$ be a Lorentzian metric. We denote the covariant components of $\metric$
with respect to a local coordinate system $x^\alpha$ by $g_{\alpha\beta}$.
The connection coefficients $\Gamma^\lambda{}_{\alpha\beta}$
and the covariant components $R_{\alpha\beta}$ of the Ricci curvature tensor
associated with $\metric$ are defined by
\bea
&&\Gamma^\lambda{}_{\alpha\beta}=
\frac{1}{2}g^{\lambda\mu}
\left(
\frac{\partial g_{\beta\mu}}{\partial x^\alpha}
+\frac{\partial g_{\alpha\mu}}{\partial x^\beta}
-\frac{\partial g_{\alpha\beta}}{\partial x^\mu}
\right),
\label{connection}
\\
&&R_{\alpha\beta}=
\frac{\partial\Gamma^\lambda{}_{\alpha\beta}}{\partial x^\lambda}
-\frac{\partial\Gamma^\lambda{}_{\beta\lambda}}{\partial x^\alpha}
+\Gamma^\lambda{}_{\alpha\beta}\Gamma^\mu{}_{\lambda\mu}
-\Gamma^\mu{}_{\alpha\lambda}\Gamma^\lambda{}_{\beta\mu}.
\label{ricci}
\eea
Here, and below, Greek indices $\alpha, \beta, \lambda,\mu$
run over the values $0,1,2,3$.
The vacuum Einstein equations may be written as
\be
R_{\alpha\beta} = 0.
\label{einsteineqco}
\ee

We look for asymptotic solutions of \eq{einsteineqco} with metrics of the form
\bea
&&\metric =
\metric\left(\frac{u(x^\alpha)}{\eps^2},\frac{y^a(x^\alpha)}{\eps},x^\alpha; \eps\right),
\label{asysol}
\\
&&\metric(\theta,\eta^a,x^\alpha;\eps) = \zero{\metric}(\theta,\eta^a,x^\alpha) + \eps \one{\metric}(\theta,\eta^a,x^\alpha)
+\eps^2 \two{\metric}(\theta,\eta^a,x^\alpha) + O(\eps^3).
\nonumber
\eea
Here, $\eps$ is a small parameter, $u$ is a phase,
$y^a$ with $a=2,3$ are spacelike transverse variables,
and
\[
\theta = \frac{u}{\eps^2},\qquad \eta^a = \frac{y^a}{\eps}
\]
are `stretched' variables. The ansatz in \eq{asysol} corresponds to a
metric that varies rapidly and strongly in the $u$-direction,
with less rapid variations in the $y^a$-directions, and slow
variations in $x^\alpha$.

We remark that the form of this ansatz is relative to a class of coordinate systems,
since an $\eps$-dependent change of coordinates can alter the
way in which the metric depends on $\eps$. It would be desirable to give a geometrically
intrinsic characterization of such an ansatz, and to carry out the expansion
in a coordinate-invariant way, but we do not attempt to do so here.

Using \eq{asysol} in the Einstein equations \eq{einsteineqco}, we find
that, to get a non-trivial solution at the leading
orders, the phase $u$ must be null to leading
order, and the transverse variables $y^a$ must be constant along the rays associated
with $u$ to
leading order. Moreover, the leading-order metric must have the form of the colliding plane
wave metric.

To carry out the expansion in detail, we use the gauge-covariance of the
Einstein equations to make a choice of coordinates that is adapted to the
metric in \eq{asysol}. We assume that $u$ is approximately null up to the order $\eps^2$,
that $y^2$, $y^3$ are constant along the rays associated
with $u$ to the order $\eps$.
Furthermore, we assume that we can extend $(u,y^2,y^3)$ to a local coordinate system
\[
(x^0,x^1,x^2,x^3)=(u,v,y^2,y^3).
\]
Then, as shown in Appendix~\ref{appA}, we can
use appropriate gauge transformations, which involve a
near identity transformation of the phase and the transverse
variables, to write a general metric \eq{asysol} whose
leading-order term has the form of the colliding plane wave metric, as
\bea
\lefteqn{
\metric = 2\zero{g}_{01}dx^0dx^1 +\zero{g}_{ab}dx^a dx^b
} \nonumber \\
&&
+\eps\left\{2\one{g}_{1a} dx^1 dx^a+\one{g}_{ab} dx^a dx^b\right\}
+\eps^2 \two{g}_{ij} dx^i dx^j
+O(\eps^3).
\label{rosen0}
\eea
In \eq{rosen0}, and below, indices $a,b,c,\dots$ take on the values $2,3$, while
indices $i,j,k,\dots$ take on the values $1,2,3$.
We raise and lower indices using the leading order metric
components; for example, we write
\[
\left(\zero{g}{}^{\alpha\beta}\right) = \left(\zero{g}{}_{\alpha\beta}\right)^{-1},\qquad
h^\alpha = \zero{g}{}^{\alpha\beta}h_\beta.
\]

We have one additional gauge freedom in \eq{rosen0}, involving a nonlinear
transformation of the phase, which we can use to set either
\[
\zero{g}_{01}=1
\]
or one of the three components
\[
\one{g}_{12},\quad \one{g}_{13},\quad \two{g}_{11}
\]
equal to zero.
We will exploit this gauge freedom later, since it is convenient in formulating a variational
principle for the asymptotic equations (see Appendix~\ref{appD}).

Using the method of multiple scales, we expand derivatives of a function $f$ with respect to $x^\mu$ as
\be
\frac{\partial}{\partial x^{\mu}}
f\left(\frac{u}{\eps^2},\frac{y^a}{\eps},x\right)=
\frac{1}{\eps^2} f_{,\theta}u_{\mu}
+\frac{1}{\eps} f_{,\bar{a}}y^a_{\mu}
+f_{,\mu},
\label{derivexp}
\ee
and then treat $(\theta,\eta^a,x)$ as independent variables.
In \eq{derivexp}, in the appendices, and below,
we use the shorthand notation
\beaa
&&f_{,\theta}
=\left.\frac{\partial f}{\partial\theta}\right|_{\eta^a, x},
\qquad
f_{,\bar{a}}
=\left.\frac{\partial f}{\partial\eta^a}\right|_{\theta, x},
\qquad
f_{,\mu}
=\left.\frac{\partial f}{\partial x^\mu}\right|_{\theta, \eta^a},\\
&&u_\mu = \frac{\partial u}{\partial x^{\mu}},
\qquad
y^a_\mu = \frac{\partial y^a}{\partial x^{\mu}}.
\eeaa

We use \eq{rosen0}--\eq{derivexp} in
\eq{connection}--\eq{ricci} and expand the result with respect to $\eps$.
We find that
\bea
&&\Gamma^\lambda{}_{\alpha\beta}=
\frac{1}{\eps^2}\mtwo{\Gamma}{}\!^\lambda{}_{\alpha\beta}
+\frac{1}{\eps}\mone{\Gamma}{}\!^\lambda{}_{\alpha\beta}
+\zero{\Gamma}{}\!^\lambda{}_{\alpha\beta}
+O(\eps),
\nonumber
\\
&&R_{\alpha\beta}=
\frac{1}{\eps^4}\mfour{R}_{\alpha\beta}
+\frac{1}{\eps^3}\mthree{R}_{\alpha\beta}
+\frac{1}{\eps^2}\mtwo{R}_{\alpha\beta}
+O(\eps^{-1}),
\label{Rexp}
\eea
where explicit expressions for the terms in the above expansions are given in Appendix~\ref{appB}.
The nonzero components of the Ricci tensor up to the order $\eps^{-2}$ are summarized in Appendix~\ref{appC}.

Using \eq{Rexp} in \eq{einsteineqco} and equating coefficients
of $\eps^{-4}$, $\eps^{-3}$ and $\eps^{-2}$ to zero,
we get
\bea
\mfour{R}_{\alpha\beta}&=&0,
\label{pert1}
\\
\mthree{R}_{\alpha\beta}&=&0,
\label{pert2}
\\
\mtwo{R}_{\alpha\beta} &=& 0.
\label{pert3}
\eea
These perturbation equations lead to a closed set of equations
for the leading-order and first-order components of the metric,
as we now explain.

The only component of \eq{pert1} that is not identically satisfied is
\be
\mfour{R}_{00} = 0.
\label{thetaconstraint}
\ee
From \eq{thetaconstraint} and \eq{Req1}, it follows that
$\zero{g}_{\alpha\beta}$ satisfies
\[
-\frac{1}{2}(\zero{g}{}\!^{ab}
\zero{g}_{ab,\theta})_{,\theta}
+\frac{1}{2}\zero{g}{}\!^{01}\zero{g}_{01,\theta}
\zero{g}{}\!^{ab} \zero{g}_{ab,\theta}
-\frac{1}{4}\zero{g}{}\!^{ac} \zero{g}_{bc,\theta}
\zero{g}{}\!^{bd} \zero{g}_{ad,\theta}=0.
\qquad\quad
\]
This equation is a constraint on the leading order metric
which involves only $\theta$-derivatives, and
has the same form as the constraint equation
for a single gravitational plane wave.

The nonzero components of the Ricci tensor at the next order in $\eps$
are
\[
\mthree{R}_{00}, \quad \mthree{R}_{0a}.
\]
The condition
\be
\mthree{R}_{0a}=0,
\label{eq0a}
\ee
yields equations for $\one{g}_{1a}$ and their derivatives with respect to $\theta$,
with nonhomogeneous terms depending on $\theta$- and $\eta$-derivatives of the leading order metric.
Thus, \eq{eq0a} provides two equations relating the first-order perturbation in the
metric to the leading order metric.

The equation
\[
\mthree{R}_{00}=0
\]
is a single equation that is homogeneous in the first-order components $\one{g}_{ab}$,
as can be seen from \eq{R00-0}. A nonzero solution of this
equation corresponds physically to a free, small-amplitude gravitational wave
of strength of the order $\eps$ propagating in the space-time
of the large-amplitude gravitational wave. Retaining a nonzero solution of
this homogenous equation would not change our final equations for the large-amplitude
wave, and for simplicity we assume that
\[
\one{g}_{ab} = 0.
\]
This assumption is also consistent with the higher-order perturbation equations.

At the next order in $\eps$, the nonzero components of the Ricci curvature
are
\[
\mtwo{R}_{00},\quad \mtwo{R}_{01},\quad \mtwo{R}_{0a},\quad \mtwo{R}_{ab}.
\]
The corresponding perturbation equations for $\mtwo{R}_{00}$ and $\mtwo{R}_{0a}$  in \eq{pert3}
give non-homogenous equations for
the second-order metric components
\[
\two{g}_{ab},\quad\two{g}_{1a},
\]
with source terms depending on
the lower-order components of the metric. These equations are satisfied
by a suitable choice of the second-order metric components, and they
are decoupled from the equations for the leading-order and first-order
metric components. We therefore do not consider them further here.

The remaining perturbation equations are
\be
\mtwo{R}_{01}=0,
\quad
\mtwo{R}_{ab}=0.
\label{apc1}
\ee
By use of the gauge freedom mentioned at the beginning of this section,
we may set
\[
\two{g}_{11}=0.
\]
In that case, \eq{apc1} provides a set of four equations relating
the leading-order and first-order metric components, and
their derivatives.

We remark that, in general, one cannot eliminate secular terms from the
asymptotic solution that are unbounded as
the `fast' phase variable $\theta$ tends to infinity.
As a result, the validity of the asymptotic equations
is restricted to a thin layer
of thickness of the order $\eps^2$ about the hypersurface $u=0$,
where $\theta=O(1)$. A global asymptotic solution can be obtained
by matching the `inner' solution inside this layer
with appropriate `outer' solutions, such as slowly varying space-times
on either side of the wave. In this
paper, however, we focus on the construction of asymptotic solutions
for localized gravitational waves, and do not consider any
matching problems.

Summarizing these results, we find that the asymptotic
solution is given by
\bea
\metric = 2\zero{g}_{01}dx^0dx^1 +\zero{g}_{ab}dx^a dx^b
+\eps\left\{2\one{g}_{1a} dx^1 dx^a\right\}
+O(\eps^2),
\label{rosen1}
\eea
where the six metric components
\[
\zero{g}_{01},\quad \zero{g}_{ab},\quad \one{g}_{1a}
\]
satisfy a system of seven equations \eq{thetaconstraint}--\eq{apc1}.
Equation \eq{thetaconstraint} is an ODE in $\theta$, and, in the case of plane-polarized
waves, we show that it is a gauge-type constraint
that is preserved by the remaining equations \eq{eq0a}--\eq{apc1}.

The explicit form
of the equations follows from the expressions in \eq{Req1}, \eq{R0a-0}--\eq{Rab-1}
for the corresponding components in the expansion of of the Ricci tensor
that appear in \eq{thetaconstraint}--\eq{apc1}.
In general, these equations are very complicated, but they can be simplified
considerably in special cases.

\subsection{Plane-polarized gravitational waves}
\label{equations}

In this section, we specialize our asymptotic solution to the case of
a plane-polarized gravitational wave that diffracts in a single direction,
and write out the resulting equations explicitly.

We choose coordinates
\be
(x^0,x^1,x^2,x^3) = (u,v,y,z)
\label{uvyz}
\ee
in which the metric has the form \eq{rosen0}. We suppose that the metric depends on the variables
$(\theta,\eta,v)$, where the phase variable $\theta$ and the transverse variable $\eta$ are defined by
\[
\theta = \frac{u}{\eps^2}, \qquad \eta = \frac{y}{\eps},
\]
and is independent of the second transverse variable $\zeta = z/\eps$.
This means that the wave diffracts only in the $y$-direction.
We could also allow the metric to depend on $z$, which would appear in the final
equations as a parameter.

We consider a leading-order metric that has the form of the colliding
plane wave metric for a plane-polarized wave, polarized
in the $(y,z)$-directions. We have seen that we need only retain
the higher order components
\[
\one{g}_{1a}, \quad \two{g}_{11}.
\]
One component of \eq{eq0a} is homogeneous in $\one{g}_{13}$,
since there is no dependence on $\zeta$, and without loss of generality
we can take
\[
\one{g}_{13}=0.
\]
We therefore take the special form of the metric \eq{rosen0} given by
\be
\metric = -2e^{-M}\left(du - \eps \Ybar dy -\frac12 \eps^2 \Tbar dv \right)dv
+e^{-U}\bigl(e^{V}dy^2+e^{-V}dz^2 \bigr) + O(\eps^2),
\label{lineernst}
\ee
where the functions $(U,V,M,\Ybar,\Tbar)$ depend on $(\theta,\eta,v)$.
We recall that we have the gauge-freedom to set either $M$, $\Ybar$,
or $\Tbar$ equal to zero.
In writing the equations, we choose to set $\Tbar=0$.

The metric \eq{lineernst} must satisfy equations \eq{thetaconstraint}--\eq{apc1}.
First, using \eq{lineernst} in \eq{thetaconstraint}, and simplifying the result, we obtain,
after some algebra, the $\theta$-constraint equation \eq{eq1}.

Next, using \eq{lineernst} and \eq{R0a-0}, we find that the only nontrivial component of \eq{eq0a}
is
\be
\mthree{R}_{02}=0.
\label{r02}
\ee
Furthermore, after the introduction of $\phi$, $\psi$ defined in \eq{defphi}--\eq{defpsi} and some algebra,
we find that \eq{r02} may be written as equation \eq{eq2}.

Finally, using \eq{lineernst} and \eq{R01-1}--\eq{Rab-1},
we find that the only nontrivial components of \eq{apc1} are
\be
\mtwo{R}_{01}=0,
\quad
\mtwo{R}_{22}=0,
\quad
\mtwo{R}_{33}=0.
\label{r2233}
\ee
After some algebra, these equations may be written as \eq{eq3}--\eq{eq5},
where the derivative $\Deta$ is defined in \eq{defdeta}

\subsection{The variational principle}

Equations \eq{thetaconstraint}--\eq{apc1}
can be derived from a variational principle, which is obtained by expanding the variational principle
for the Einstein equations. In particular, the variational principle provides a check on the
algebra used to derive the equations.
(We have also checked the final results by use of MAPLE.)

In order to formulate a variational principle,
it is necessary to retain the component
\[
\two{g}_{11}.
\]
Variations with respect to this component
yield the gauge-type constraint \eq{thetaconstraint}, and it may be set to zero after taking variations
with respect to it.

The general form of the asymptotic variational principle is given in Appendix~\ref{appD}.
Here, we specialize it to the case of a plane-polarized gravitational wave
considered in Section~\ref{equations}. Using the metric \eq{lineernst} in \eq{I-2}
we find, after some algebra, that the variational principle \eq{varprin} for
\eq{eq1}--\eq{eq5} is
\[
\delta \Sone = 0,\qquad \Sone = \int \Lmtwo\,d\theta\, dv \,d\eta,
\]
where the Lagrangian $\Lmtwo$ may be written as
\bea
&&
\Lmtwo
=
e^{-U}\Big\{
2 M_{\theta v} +4 U_{\theta v}-V_{\theta} V_{v} -3 U_{\theta} U_{v}
\nonumber \\
&&
\qquad
-\Tbar_{\theta\theta}
+\Tbar_{\theta}(M_{\theta} + 2 U_{\theta})
+\Tbar(M_{\theta\theta}+2 U_{\theta\theta}
-{\textstyle\frac32} U_{\theta}^{2} -{\textstyle\frac12} V_{\theta}^{2})
\nonumber \\
&&
\qquad
+ e^{-(U+V+M)}\left[
-2\Deta\phi - \Deta\psi + {\textstyle\frac32} \phi^2
+2\phi\psi + \psi^2\right]
\Big\} .
\label{I-1}
\eea
Variations of $\Sone$ with respect
to the second-order metric component $\Tbar$
lead to the $\theta$-constraint equation \eq{eq1}.
Variations with respect to the first-order metric component $\Ybar$,
which appears in $\phi$, $\psi$, and $D_\eta$,
lead to equation \eq{eq2}.
Variations with respect to $U$, $V$, $M$
lead to the evolution equations \eq{eq3}--\eq{eq5},
after we set $\Tbar=0$.

\section{Properties of the equations}
\label{properties}
\setcounter{equation}{0}

In this section, we study some properties of the equations for plane-polarized
waves that we have derived in the previous section. We write out the structure of the highest-order
derivatives that appear in the equations, and use this
to formulate a reasonable IBVP for them. We also show that the
$\theta$-constraint equation is preserved by the evolution in $v$,
and that the linearized equations are consistent with the
linearized equations for gravitational waves in the parabolic
approximation.

\subsection{Structure of equations and an IBVP}
\label{ibvp}

In this subsection, we consider the structure of equations \eq{eq1}--\eq{eq5} in more detail.
The first equation, \eq{eq1}, is an ODE with respect to $\theta$ relating $(U,V,M)$.
As we show in the next section,
this equation is a gauge-type constraint which holds for all $v$
if it holds for $v=0$, say. We may therefore
neglect this equation provided that the
initial data at $v=0$ is compatible with it.

The remaining equations \eq{eq2}--\eq{eq5} form a system of equations for $(U,V,M,\Ybar)$. In
order to exhibit their structure, we rewrite them in a way that shows explicitly how
the highest, second-order, derivatives appear.

Using \eq{defdeta}--\eq{defpsi}, we may rewrite equation \eq{eq2} as
\bea
&&\Ybar_{\theta\theta} - \left\{(V+M)_\theta\Ybar\right\}_\theta
+\left\{\frac{1}{2} U_\theta^2 - U_\theta V_\theta -\frac{1}{2} V_\theta^2\right\} Y
\nonumber
\\
&&\qquad\qquad = \left(U+V+M\right)_{\theta\eta} + M_\eta U_\theta - (U+V)_\eta V_\theta.
\label{Yode}
\eea
If $(U,V,M)$ are assumed known, then this equation is a linear ODE in $\theta$ for $\Ybar$.

We may rewrite equations \eq{eq3}--\eq{eq5} as
\bea
&&(U+V+M)_{\theta v} - \frac{1}{2}\left(U + V\right)_\theta\left(U + V\right)_v
\nonumber\\
&&\qquad=\frac{1}{2} e^{-(U+V+M)} \left\{-\frac{1}{2}\phi^2 - \phi\psi\right\},
\label{ee1}\\
&&(U+V)_{\theta v} - \frac{1}{2}U_\theta \left(U + V\right)_v - \frac{1}{2}U_v\left(U + V\right)_\theta,
\nonumber\\
&&\qquad = \frac{1}{2} e^{-(U+V+M)} \left\{\Deta^2(U+V) - \phi\psi - \psi^2\right\}
\label{ee2}\\
&&U_{\theta v} - U_\theta U_v
\nonumber \\
&&\qquad = \frac{1}{2}e^{-(U+V+M)}\Bigl\{\Deta^2(U+V+M) - \Deta\left(e^{U}\Ybar_\theta\right)
\nonumber\\
&&\qquad\qquad\qquad\qquad\qquad - \frac{1}{2} \phi^2 - \phi\psi - \psi^2\Bigr\}.
\label{ee3}
\eea
Examining the terms that involve second order derivatives, we see that these equations
consist of a $(1+2)$-dimensional wave equation in $(\theta,\eta,v)$ for $(U+V)$,
and two $(1+1)$-dimensional wave equations in $(\theta,v)$, for $(U+V+M)$ and $U$,
in which $(\theta,v)$ are characteristic coordinates.
An additional second-order derivative term, proportional to $\Deta \Ybar_\theta$,
appears in the equation for $U$. The function $\Ybar$ is also coupled with
the evolution equations through the dependence of the transverse derivative $\Deta$, given in \eq{defdeta},
on $\Ybar$.

In order to specify a unique solution of these equations, we expect that we
need to supplement the ODE \eq{Yode} with data for $\Ybar$ and $\Ybar_\theta$ on $\theta=0$, say,
and the evolution equations \eq{ee1}--\eq{ee3} with characteristic initial data for $(U,V,M)$
on $\theta = 0$ and $v=0$.
Thus, a reasonable IBVP for \eq{Yode}--\eq{ee3} in the region $\theta > 0$,
$-\infty<\eta < \infty$, and $v > 0$ is
\beaa
&&(U,V,M) = (U_0,V_0,M_0) \qquad\mbox{on $v=0$},\\
&&(U,V,M) = (U_1,V_1,M_1) \qquad\mbox{on $\theta=0$},\\
&&(Y,Y_\theta) = (Y_0,Y_1)\qquad\qquad\quad\mbox{on $\theta=0$}.
\eeaa
Here, $(U_0,V_0,M_0)$ are given functions of $(\theta,\eta)$ that satisfy the
constraint
\[
U_{0\theta\theta} - \frac{1}{2}\left(U_{0\theta}^2 + V_{0\theta}^2\right) + U_{0\theta}M_{0\theta} = 0,
\]
and $(U_1,V_1,M_1, Y_0,Y_1)$ are given functions of $(\eta,v)$. This data
may be interpreted as initial data for the state of the wave
on the hypersurface $v=0$, and boundary data on the leading wavefront $\theta=0$.

\subsection{The constraint equations}
\label{constraints}

In this section, we show that the constraint equation \eq{eq1} is preserved.

\begin{proposition}
\label{prop:con}
Suppose that $(U,V,M, Y)$ are smooth functions that satisfy
\eq{eq2}--\eq{eq5}. Let
\be
F = U_{\theta\theta} - \frac{1}{2}\left(U_\theta^2 + V_\theta^2\right) + U_\theta M_\theta.
\label{eeF}
\ee
Then
\be
F_v = U_v F.
\label{conF}
\ee
\end{proposition}

\bigskip\noindent
\textbf{Proof.}
We write the evolution equations \eq{eq3}--\eq{eq5} as
\bea
&&
U_{\theta v}- U_{\theta}U_{v}
=\frac12 e^{-(U+V+M)}A,
\label{eeA}\\
&&
V_{\theta v}
-\frac{1}{2}\left( U_{\theta}V_{v}+U_{v}V_{\theta}\right)
=\frac12 e^{-(U+V+M)}B,
\label{eeB}\\
&&
M_{\theta v}
+\frac{1}{2}\left( U_{\theta}U_{v}-V_{v}V_{\theta}\right)
=\frac12 e^{-(U+V+M)}C,
\label{eeC}
\eea
where, using the notation defined in \eq{defdeta}--\eq{defpsi},
\bea
&&
A=\Deta\phi+\Deta\psi-\frac12 \phi^2- \phi\psi -\psi^2,
\label{defA}
\\
&&
B=-\Deta\phi+\frac12 \phi^2,
\label{defB}
\\
&&
C=-\Deta\psi -\frac12 \phi^2 +\psi^2.
\label{defC}
\eea

Differentiating \eq{eeF} with respect to $v$, and using \eq{eeA}--\eq{eeC} and \eq{eeF}
to replace $U_{\theta v}$, $V_{\theta v}, M_{\theta v}$ and $U_{\theta\theta}$ in the
result, we find that
\be
F_v = U_v F + \frac{1}{2}e^{-(U+V+M)} D,
\label{Feq}
\ee
where
\be
D = A_\theta - (U+V)_\theta A - V_\theta B + U_\theta C.
\label{defD}
\ee
Differentiating equation \eq{defA} for $A$ with respect to $\theta$,
introducing the commutator  $\left[\partial_\theta,\Deta\right]$
of $\partial_\theta$ and $D_\eta$, and using equations \eq{defpsi} and
\eq{eq2},
we compute that
\bea
A_\theta &=& \partial_\theta\left\{\Deta(\phi + \psi) - \frac{1}{2}\phi^2 - \phi\psi -\psi^2\right\}
\nonumber\\
&=& \Deta (\phi+\psi)_\theta  + \left[\partial_\theta,\Deta\right](\phi+\psi)
- (\phi+\psi)\phi_\theta - (\phi+2\psi)\psi_\theta
\nonumber\\
&=& \Deta \left\{\psi(U+V)_\theta\right\}
\nonumber\\
&&\qquad - (\phi+\psi)\phi_\theta - (\phi+2\psi)\psi_\theta + \left[\partial_\theta,\Deta\right](\phi+\psi)
\nonumber\\
&=& (U+V)_\theta \Deta\psi +  \psi\Deta(U+V)_\theta
\nonumber\\
&&\qquad
- (\phi+\psi)\phi_\theta - (\phi+2\psi)\psi_\theta + \left[\partial_\theta,\Deta\right](\phi+\psi)
\nonumber\\
&=& (U+V)_\theta \Deta\psi + \psi\partial_\theta \Deta (U+V)
- \psi\left[\partial_\theta,\Deta\right](U+V)
\nonumber\\
&&\qquad
- (\phi+\psi)\phi_\theta - (\phi+2\psi)\psi_\theta + \left[\partial_\theta,\Deta\right](\phi+\psi)
\nonumber\\
&=&
(U+V)_\theta \Deta\psi + \psi\psi_\theta - (\phi+\psi)\phi_\theta - (\phi+2\psi)\psi_\theta
\nonumber\\
&&\qquad
+ \left[\partial_\theta,\Deta\right](\phi+\psi)
- \psi\left[\partial_\theta,\Deta\right](U+V)
\nonumber\\
&=& (U+V)_\theta \Deta\psi - (\phi+\psi)(\phi+\psi)_\theta
\nonumber\\
&&\qquad
+ \left[\partial_\theta,\Deta\right](\phi+\psi)
- \psi\left[\partial_\theta,\Deta\right](U+V)
\nonumber\\
&=& (U+V)_\theta\left\{\Deta\psi - \psi(\phi+\psi)\right\}
\nonumber\\
&&\qquad
+ \left[\partial_\theta,\Deta\right](\phi+\psi)
- \psi\left[\partial_\theta,\Deta\right](U+V).
\label{eqAtheta}
\eea
From the definition of $\Deta$ in \eq{defdeta}, it follows that
\[
\left[\partial_\theta,\Deta\right] = U_\theta \Deta + e^{U}\Ybar_\theta \partial_\theta.
\]
Hence, using \eq{defphi}--\eq{defpsi} and \eq{eq2}, we get
\beaa
\left[\partial_\theta,\Deta\right](\phi+\psi)
&=& U_\theta \Deta(\phi+\psi) + e^U \Ybar_\theta(\phi+\psi)_\theta\\
&=& U_\theta \Deta(\phi+\psi) + \psi e^{U}\Ybar_\theta(U+V)_\theta,\\
\left[\partial_\theta,\Deta\right](U+V)
&=& U_\theta \Deta(U+V) + e^U Y_\theta(U+V)_\theta\\
&=& U_\theta \psi + e^{U}\Ybar_\theta(U+V)_\theta.
\eeaa
Using these equations in \eq{eqAtheta}, and simplifying the result, we
find that
\[
A_\theta = (U+V)_\theta\left(\Deta \psi - \phi\psi - \psi^2\right)
+ U_\theta\left(\Deta\phi + \Deta\psi - \psi^2\right).
\]
Finally, using this equation and \eq{defA}--\eq{defC} in \eq{defD},
and simplifying the result, we find
that $D=0$. It follows from \eq{Feq} that $F$ satisfies \eq{conF}.\quad$\Box$

\subsection{Linearization}
\label{linearization}

We consider the small-amplitude limit of \eq{eq1}--\eq{eq5} in which
\[
U,V,M,\Ybar \to 0.
\]
From the constraint equation \eq{eq1}, we have $U = O(V^2)$, so $U$ is of higher
order in a linearized approximation and can be neglected completely.
From \eq{defdeta}--\eq{defphi}, we also have in this approximation that
\beaa
&&\Deta = \partial_\eta,\\
&&\phi =  M_\eta - \Ybar_\theta,\\
&&\psi = V_\eta.
\eeaa
Linearization of \eq{eq2} yields
\[
(\phi + \psi)_\theta = 0,
\]
while linearization the evolution equations \eq{eq3}--\eq{eq5} gives
\beaa
&&0 = \frac{1}{2}(\phi+\psi)_\eta,\\
&&V_{\theta v} = -\frac{1}{2}\phi_\eta,\\
&&M_{\theta v} = -\frac{1}{2} \psi_{\eta}.
\eeaa

Neglecting some arbitrary functions of integration for simplicity,
we find that these
equations are satisfied if $M=-V$ and $\Ybar=0$.
In that case, $V$ satisfies an equation of the form \eq{pareq}:
\[
V_{\theta v} = \frac{1}{2} V_{\eta\eta}.
\]
The corresponding linearized metric is given by
\[
\metric = 2(1-V)du dv + (1+V)dy^2 + (1-V) dz^2.
\]
Linearization of the Einstein equations, with a suitable choice of gauge, leads to a set of
linear wave equations for the metric components. One can verify that this linearization
of our asymptotic solution agrees with what is obtained by an application of the parabolic approximation,
described in Section~\ref{sec:wave}, to the linearized Einstein equations.

\appendix

\section{Gauge transformations}
\setcounter{equation}{0}
\label{appA}

We consider a Lorentzian metric $\metric$, and local coordinates $x^\alpha$, in which
\[
\metric=g_{\alpha\beta}dx^\alpha\,dx^\beta.
\]
We denote the contravariant form of $\metric$ by $\metric^\sharp$.

We suppose that the metric $\metric$ depends upon a small parameter
$\eps$, and there exist independent functions $u(x)$, $y^a(x)$, $a=2,3$, such that
as $\eps\to 0$ we have
\bea
&&
\metric^\sharp(du,du)=O(\eps^2),
\label{cond0}
\\ &&
\metric^\sharp(du,dy^a)=O(\eps), \quad a=2,3.
\label{cond2}
\eea
The first condition states that $du$ is an approximate null form up to the order $\eps^2$.
The second condition states that
$y^a$ is approximately constant along the rays associated with $u$.

In this appendix, we reduce such a metric to the standard asymptotic form,
given in \eq{rosen0}. It is convenient to consider two different
types of coordinate transformations.
First, we show that it is possible to make a near-identity
change of coordinates, to a local
coordinate system $x^\alpha$ with
\[
x^0=u+O(\eps^3),\qquad x^a=y^a+O(\eps^2),
\]
such that
\bea
\lefteqn{
\metric = 2\zero{g}_{01}dx^0dx^1 +\zero{g}_{ab}dx^a dx^b
} \nonumber \\
&&
+\eps\left\{2\one{g}_{1a} dx^1 dx^a+\one{g}_{ab} dx^a dx^b\right\}
+\eps^2 \two{g}_{ij} dx^i dx^j
+O(\eps^3).
\label{rosen}
\eea
Second, we show that the form of the metric can be further simplified by
use of a coordinate $x^0$ such that
\[
x^0=u+O(\eps^2).
\]
As before, indices $\alpha,\beta,\dots$ take on the values $0,1,2,3$; indices $a,b,c,\dots$
take on the values $2,3$; and indices $i,j,k,\dots$ take on the values
$1,2,3$.

Let $x^\alpha$ be any local coordinate system in which
\[
x^0 = u, \qquad x^a = y^a.
\]
Then the conditions \eq{cond0}, \eq{cond2} are equivalent to
\be
\zero{g}\!^{00}=\one{g}\!^{00}=0,
\quad
\zero{g}\!^{0a}=0.
\label{cond}
\ee
These conditions imply immediately that
\be
\zero{g}_{11}=0, \qquad \zero{g}_{1a}=0.
\label{rest}
\ee
Moreover, we have
\be
\det (g_{ab}) > 0,
\label{det}
\ee
which says that $y^a$ are space-like independent coordinates.

The expansion \eq{asysol}, with the additional restriction \eq{rest}, is invariant under
a coordinate transformation of the form
\be
x^\alpha
\rightarrow
x^\alpha+\eps^2 \Psi^\alpha\left(\frac{x^0}{\eps^2},\frac{x^a}{\eps},x;\eps\right),
\label{tr}
\ee
\beaa
&&
\Psi^0(\theta,\eta^a,x;\eps)
= \eps \one{\Psi}{}\!^0\left(\theta,\eta^a,x\right)
+ \eps^2 \two{\Psi}{}\!^0\left(\theta,\eta^a,x\right)
+O(\eps^3),
\\ &&
\Psi^i(\theta,\eta^a,x;\eps)
= \zero{\Psi}{}\!^i\left(\theta,\eta^a,x\right)
+ \eps \one{\Psi}{}\!^i\left(\theta,\eta^a,x\right)
+ \eps^2 \two{\Psi}{}\!^i\left(\theta,\eta^a,x\right)
+O(\eps^3),
\eeaa
where, as usual,
\[
\theta = \frac{x^0}{\eps},\qquad \eta^a = \frac{x^a}{\eps}.
\]

To begin with, we perform a change of coordinates of the form
\be
x^0 \rightarrow x^0,
\qquad
x^i \rightarrow x^i+\eps^2\zero{\Psi}{}\!^i(\theta,\eta^c,x).
\label{tt0}
\ee
At order zero in $\eps$, the metric components transform according to
\beaa
\zero{g}_{00}
& \rightarrow &
\zero{g}_{00}
+2 \zero{\Psi}{}\!^k_{,\theta}\zero{g}_{0k}
+ \zero{\Psi}{}\!^k_{,\theta} \zero{\Psi}{}\!^l_{,\theta}\zero{g}_{kl},
\\
\zero{g}_{0i}
& \rightarrow &
\zero{g}_{0i}
+\zero{\Psi}{}\!^k_{,\theta}\zero{g}_{ik},
\\
\zero{g}_{ij}
& \rightarrow &
\zero{g}_{ij}.
\eeaa
Recalling \eq{det}, the matrix $(g_{ab})$
is invertible, and we can choose \[\zero{\Psi}{}\!^i_{,\theta}\] so that
\be
\zero{g}_{00}=0, \qquad \zero{g}_{0a}=0.
\ee
Thus, the zero order metric $\zero{\metric}$ takes the form given in \eq{rosen}.

Next, we consider the effect of the coordinate
transformation \eq{tr} on the
metric components at order $\eps$.
Under the action of the transformation
\be
x^\alpha \rightarrow x^\alpha+\eps^3\one{\Psi}{}\!^\alpha(\theta,\eta^c,x),
\ee
the form \eq{rosen}
of the metric is unchanged at order zero.
At order one, the components transforms according to
\beaa
\one{g}_{00}
& \rightarrow &
\one{g}_{00}
+2 \one{\Psi}{}\!^1_{,\theta}\zero{g}_{01},
\\
\one{g}_{01}
& \rightarrow &
\one{g}_{01}
+\one{\Psi}{}\!^0_{,\theta}\zero{g}_{01},
\\
\one{g}_{0a}
& \rightarrow &
\one{g}_{0a}
+\one{\Psi}{}\!^b_{,\theta}\zero{g}_{ab},
\\
\one{g}_{ij}
& \rightarrow &
\one{g}_{ij}.
\eeaa
Choosing appropriately $\one{\Psi}{}\!^\alpha$,
these transformations can be used to make
\be
\one{g}_{0\alpha}=0.
\ee
Moreover, the condition \eq{cond} implies
\be
\one{g}_{11} \equiv -\left(\zero{g}_{01}\right)^2 \one{g}\!^{00}=0.
\ee

Next, we simplify the metric components at order $\eps^2$.
Under the action of the transformation
\be
x^\alpha \rightarrow
x^\alpha+\eps^4\two{\Psi}{}\!^\alpha(\theta,\eta^c,x),
\label{trans1}
\ee
the form of the metric is unchanged at order zero and one.
At order two, the metric components transform according to
\beaa
\two{g}_{00}
& \rightarrow &
\two{g}_{00}
+2 \two{\Psi}{}\!^1_{,\theta}\zero{g}_{01},
\\
\two{g}_{01}
& \rightarrow &
\two{g}_{01}
+\two{\Psi}{}\!^0_{,\theta}\zero{g}_{01},
\\
\two{g}_{0a}
& \rightarrow &
\two{g}_{0a}
+\two{\Psi}{}\!^b_{,\theta}\zero{g}_{ab},
\\
\two{g}_{ij}
& \rightarrow &
\two{g}_{ij}.
\eeaa
Choosing $\two{\Psi}{}\!^\alpha$ appropriately,
these transformations can be used to make
\be
\two{g}_{0\alpha}=0,
\ee
and the metric can be written as in \eq{rosen}.

If we allow for changes of order $\eps^2$ in $x^0$, then
we still have the freedom to use a transformation of the form
\be
x^0 \rightarrow \eps^2 \Psi(\theta,\eta^c,x).
\label{trans2}
\ee
The substantial difference between \eq{trans1} and \eq{trans2} is that,
unlike the former one, the latter involves a nonlinear transformation
of the phase $\theta$. This transformation can be interpreted as a transformation to
characteristic coordinates in the space $(\theta,\eta^a,x)$.
The gauge \eq{trans2} acts only on the following components:
\[
\zero{g}_{01} \rightarrow
\Psi_{,\theta}\zero{g}_{01},
\qquad
\one{g}_{1a} \rightarrow
\one{g}_{1a}+\Psi_{,\bar{a}}\zero{g}_{01},
\qquad
\two{g}_{11} \rightarrow
\two{g}_{11}+2 \Psi_{,\theta}\zero{g}_{01}.
\]
We can use \eq{trans2} to set $\zero{g}_{01}=1$ or one of the three
components
\[
\one{g}_{12},\quad \one{g}_{13},\quad \two{g}_{11}
\]
equal to zero.
More generally, it is possible to choose $\Psi$ so that an appropriate relationship involving
all of these components is satisfied.

The contravariant form of
the metric tensor in \eq{rosen} is
\bea
\lefteqn{
\metric^\sharp
=2\zero{g}{}\!^{01}\partial_0\partial_1
+\zero{g}{}\!^{ab}\partial_a\partial_b
-\eps\left\{
2\one{g}{}\!^{0a}\partial_0\partial_a
+\one{g}{}\!^{ab}\partial_a\partial_b\right\} \nonumber
} \nonumber
\\
&&
+\eps^2\left\{
\two{h}{}\!^{00}\partial_0^2
+2\two{h}{}\!^{0a}\partial_0\partial_a
+\two{h}{}\!^{ab}\partial_a\partial_b\right\}
+O(\eps^3),
\label{contra}
\eea
with
\beaa
&&
\two{h}{}\!^{00} = -(\two{g}{}\!^{00}-\zero{g}_{cd}\one{g}{}\!^{0c}\one{g}{}\!^{0d}),
\\ &&
\two{h}{}\!^{0a} = -(\two{g}{}\!^{0a}-\zero{g}_{cd}\one{g}{}\!^{0c}\one{g}{}\!^{ad}),
\\ &&
\two{h}{}\!^{ab} = -(\two{g}{}\!^{ab}-\zero{g}_{cd}\one{g}{}\!^{ac}\one{g}{}\!^{bd}).
\eeaa
In \eq{contra}, $\zero{g}{}\!^{01}$ and $\zero{g}{}\!^{ab}$ are defined by
\[
\zero{g}{}\!^{01}\zero{g}_{01}=1,
\qquad
\zero{g}{}\!^{ac}\zero{g}_{cb}=\delta^a_b.
\]
Moreover, $\zero{g}{}\!^{01}$ is used to raise the sub-index $1$, and
$\zero{g}{}\!^{ab}$ is used to raise the sub-indices $2,3$.
Thus, for example, we have
\[
\one{g}{}\!^{00}=(\zero{g}{}\!^{01})^2\one{g}_{11},
\qquad
\two{g}{}\!^{0a}=\zero{g}{}\!^{01}\zero{g}{}\!^{ab}\two{g}_{1b}.
\]
It follows from \eq{contra} that the transformed phase $u=x^0$ and
transverse variables $y^a=x^a$ satisfy \eq{cond0}--\eq{cond2}.

\section{Expansion of connection coefficients and Ricci tensor}
\setcounter{equation}{0}
\label{appB}

We look for a multiple-scale expansion as $\eps \to 0$
of the metric components of the form
\bea
g_{\alpha\beta} &=&
g_{\alpha\beta}\left(\frac{u(x)}{\eps^2},\frac{y^a(x)}{\eps},x;\eps\right),
\label{exp1}
\\
g_{\alpha\beta}\left(\theta,\eta^a,x;\eps\right)
&=& \zero{g}_{\alpha\beta}\left(\theta,\eta^a,x\right)
+\eps \one{g}_{\alpha\beta}\left(\theta,\eta^a,x\right)
+\eps^2 \two{g}_{\alpha\beta}\left(\theta,\eta^a,x\right)
+O(\eps^3).
\nonumber
\eea
The contravariant metric components $g^{\alpha\beta}$ satisfy
\[
g^{\alpha\mu}g_{\mu\beta}=\delta^\alpha_\beta.
\]
Expansion of this equation in a power series in $\eps$
gives
\bea
g^{\alpha\beta}
& =& \zero{h}\!^{\alpha\beta}
+\eps \one{h}\!^{\alpha\beta}
+\eps^2 \two{h}\!^{\alpha\beta} + O(\eps^3)
\nonumber
\\
& := & \zero{g}{}\!^{\alpha\beta}
-\eps \one{g}{}\!^{\alpha\beta}
-\eps^2\left(\two{g}{}\!^{\alpha\beta}
-\zero{g}_{\mu\nu}\one{g}{}\!^{\alpha\mu}\one{g}{}\!^{\beta\nu}\right)
+O(\eps^3),
\label{exp2}
\eea
where $\zero{g}{}\!^{\alpha\beta}$ is the inverse
of $\zero{g}_{\alpha\beta}$.
In \eq{exp2} and below, we use the leading order metric
components to raise indices, so that
\be
\one{g}{}\!^{\alpha\beta}=
\zero{g}{}\!^{\alpha\mu}\zero{g}{}\!^{\beta\nu}
\one{g}_{\mu\nu},
\qquad
\two{g}{}\!^{\alpha\beta}=
\zero{g}{}\!^{\alpha\mu}\zero{g}{}\!^{\beta\nu}
\two{g}_{\mu\nu}.
\label{g1}
\ee
With this notation, note that, for instance,
the first order term in the expansion of the covariant
metric components $g^{\alpha\beta}$ with respect to $\eps$ is not $\one{g}{}\!^{\alpha\beta}$
but $\one{h}\!^{\alpha\beta}=-\one{g}{}\!^{\alpha\beta}$.

The expansion of the connection coefficients (or Christoffel symbols) is
\beaa
&&
\Gamma^\lambda{}_{\alpha\beta}=
\frac{1}{\eps^2}\mtwo{\Gamma}{}\!^\lambda{}_{\alpha\beta}
+\frac{1}{\eps}\mone{\Gamma}{}\!^\lambda{}_{\alpha\beta}
+\zero{\Gamma}{}\!^\lambda{}_{\alpha\beta}
+O(\eps),
\eeaa
where
\bea
\mtwo{\Gamma}{}\!^\lambda{}_{\alpha\beta}
&=&
\frac{1}{2}\zero{g}{}\!^{\lambda\mu}
\left(
\zero{g}_{\beta\mu,\theta}u_\alpha
+\zero{g}_{\alpha\mu,\theta}u_\beta
-\zero{g}_{\alpha\beta,\theta}u_\mu
\right),
\nonumber
\\
\mone{\Gamma}{}\!^\lambda{}_{\alpha\beta}
&=&
\frac{1}{2}\zero{g}{}\!^{\lambda\mu}
\left(
\zero{g}_{\beta\mu,\bar{a}}y^a_{\alpha}
+\zero{g}_{\alpha\mu,\bar{a}}y^a_{\beta}
-\zero{g}_{\alpha\beta,\bar{a}}y^a_{\mu}
\right)
\nonumber
\\
&&
+\frac{1}{2}\zero{g}{}\!^{\lambda\mu}
\left(
\one{g}_{\beta\mu,\theta}u_\alpha
+\one{g}_{\alpha\mu,\theta}u_\beta
-\one{g}_{\alpha\beta,\theta}u_\mu
\right)
\nonumber
\\
&&
-\frac{1}{2}\one{g}{}\!^{\lambda\mu}
\left(
\zero{g}_{\beta\mu,\theta}u_\alpha
+\zero{g}_{\alpha\mu,\theta}u_\beta
-\zero{g}_{\alpha\beta,\theta}u_\mu
\right),
\nonumber
\\
\zero{\Gamma}{}\!^\lambda{}_{\alpha\beta}
&=&
\frac{1}{2}\zero{g}{}\!^{\lambda\mu}
\left(
\zero{g}_{\beta\mu,\alpha}
+\zero{g}_{\alpha\mu,\beta}
-\zero{g}_{\alpha\beta,\mu}
\right)
\nonumber
\\
&&
+\frac{1}{2}\zero{g}{}\!^{\lambda\mu}
\left(
\one{g}_{\beta\mu,\bar{a}}y^a_\alpha
+\one{g}_{\alpha\mu,\bar{a}}y^a_\beta
-\one{g}_{\alpha\beta,\bar{a}}y^a_\mu
\right)
\nonumber
\\
&&
-\frac{1}{2}\one{g}{}\!^{\lambda\mu}
\left(
\zero{g}_{\beta\mu,\bar{a}}y^a_\alpha
+\zero{g}_{\alpha\mu,\bar{a}}y^a_\beta
-\zero{g}_{\alpha\beta,\bar{a}}y^a_\mu
\right)
\nonumber
\\
&&
+\frac{1}{2}\zero{g}{}\!^{\lambda\mu}
\left(
\two{g}_{\beta\mu,\theta}u_\alpha
+\two{g}_{\alpha\mu,\theta}u_\beta
-\two{g}_{\alpha\beta,\theta}u_\mu
\right)
\nonumber
\\
&&
-\frac{1}{2}\one{g}{}\!^{\lambda\mu}
\left(
\one{g}_{\beta\mu,\theta}u_\alpha
+\one{g}_{\alpha\mu,\theta}u_\beta
-\one{g}_{\alpha\beta,\theta}u_\mu
\right)
\nonumber
\\
&&
+\frac{1}{2} \two{h}\!^{\lambda\mu}
\left(
\zero{g}_{\beta\mu,\theta}u_\alpha
+\zero{g}_{\alpha\mu,\theta}u_\beta
-\zero{g}_{\alpha\beta,\theta}u_\mu
\right).
\nonumber
\eea

The expansion of the Ricci tensor is
\beaa
&&
R_{\alpha\beta}=
\frac{1}{\eps^4}\mfour{R}_{\alpha\beta}
+\frac{1}{\eps^3}\mthree{R}_{\alpha\beta}
+\frac{1}{\eps^2}\mtwo{R}_{\alpha\beta}
+O(\eps^{-1}),
\eeaa
where
\bea
\mfour{R}_{\alpha\beta}
&=&
\mtwo{\Gamma}{}\!^\mu{}_{\alpha\beta,\theta}u_\mu
-\mtwo{\Gamma}{}\!^\mu{}_{\beta\mu,\theta}u_\alpha
+\mtwo{\Gamma}{}\!^\mu{}_{\alpha\beta}
\mtwo{\Gamma}{}\!^\nu{}_{\mu\nu}
-\mtwo{\Gamma}{}\!^\mu{}_{\alpha\nu}
\mtwo{\Gamma}{}\!^\nu{}_{\beta\mu},
\nonumber
\\
\mthree{R}_{\alpha\beta}
&=&
\mone{\Gamma}{}\!^\mu{}_{\alpha\beta,\theta}u_\mu
-\mone{\Gamma}{}\!^\mu{}_{\beta\mu,\theta}u_\alpha
+\mtwo{\Gamma}{}\!^\mu{}_{\alpha\beta,\bar{a}}y^a_{\mu}
-\mtwo{\Gamma}{}\!^\mu{}_{\beta\mu,\bar{a}}y^a_{\alpha}
\nonumber
\\
&&
+\mtwo{\Gamma}{}\!^\mu{}_{\alpha\beta}
\mone{\Gamma}{}\!^\nu{}_{\mu\nu}
+\mone{\Gamma}{}\!^\mu{}_{\alpha\beta}
\mtwo{\Gamma}{}\!^\nu{}_{\mu\nu}
-\mtwo{\Gamma}{}\!^\mu{}_{\alpha\nu}
\mone{\Gamma}{}\!^\nu{}_{\beta\mu}
-\mone{\Gamma}{}\!^\mu{}_{\alpha\nu}
\mtwo{\Gamma}{}\!^\nu{}_{\beta\mu},
\nonumber
\\
\mtwo{R}_{\alpha\beta}
&=&
\zero{\Gamma}{}\!^\mu{}_{\alpha\beta,\theta}u_\mu
-\zero{\Gamma}{}\!^\mu{}_{\beta\mu,\theta}u_\alpha
+\mone{\Gamma}{}\!^\mu{}_{\alpha\beta,\bar{a}}y^a_{\mu}
-\mone{\Gamma}{}\!^\mu{}_{\beta\mu,\bar{a}}y^a_{\alpha}
\nonumber
\\
&&
+\mtwo{\Gamma}{}\!^\mu{}_{\alpha\beta,\mu}
-\mtwo{\Gamma}{}\!^\mu{}_{\beta\mu,\alpha}
+\mtwo{\Gamma}{}\!^\mu{}_{\alpha\beta}
\zero{\Gamma}{}\!^\nu{}_{\mu\nu}
+\zero{\Gamma}{}\!^\mu{}_{\alpha\beta}
\mtwo{\Gamma}{}\!^\nu{}_{\mu\nu}
\nonumber
\\
&&
-\mtwo{\Gamma}{}\!^\mu{}_{\alpha\nu}
\zero{\Gamma}{}\!^\nu{}_{\beta\mu}
-\zero{\Gamma}{}\!^\mu{}_{\alpha\nu}
\mtwo{\Gamma}{}\!^\nu{}_{\beta\mu}
+\mone{\Gamma}{}\!^\mu{}_{\alpha\beta}
\mone{\Gamma}{}\!^\nu{}_{\mu\nu}
-\mone{\Gamma}{}\!^\mu{}_{\alpha\nu}
\mone{\Gamma}{}\!^\nu{}_{\beta\mu}.
\nonumber
\eea

\section{Nonzero connection coefficients and Ricci tensor components}
\setcounter{equation}{0}
\label{appC}

It this appendix, we write out expressions for the nonzero connection
coefficients and Ricci tensor components for a metric with the form
given in \eq{rosen0}.

Nonzero connection coefficients at order $\eps^{-2}$:
\beaa
&&
\mtwo{\Gamma}{}\!^0{}_{00}
=\zero{g}{}\!^{01}\zero{g}_{01,\theta},
\quad
\mtwo{\Gamma}{}\!^{1}{}_{ab}
=-\frac{1}{2}\zero{g}{}\!^{01}\zero{g}_{ab,\theta},
\quad
\mtwo{\Gamma}{}\!^{a}{}_{0b}
=\frac{1}{2}\zero{g}{}\!^{ac}\zero{g}_{bc,\theta}.
\eeaa
Nonzero connection coefficients at order $\eps^{-1}$:
\beaa
&&
\mone{\Gamma}{}\!^{0}{}_{0a}
=\frac{1}{2}\zero{g}{}\!^{01}
(\zero{g}_{01,\bar{a}}+\one{g}_{1a,\theta})
-\frac{1}{2}\one{g}{}\!^{0b}\zero{g}_{ab,\theta},
\\
&&
\mone{\Gamma}{}\!^{1}{}_{1a}
=\frac{1}{2}\zero{g}{}\!^{01}
(\zero{g}_{01,\bar{a}}-\one{g}_{1a,\theta}),
\quad
\mone{\Gamma}{}\!^{1}{}_{ab}
=-\frac{1}{2}\zero{g}{}\!^{01}
\one{g}_{ab,\theta},
\\
&&
\mone{\Gamma}{}\!^{a}{}_{01}
=-\frac{1}{2}\zero{g}{}\!^{ac}
(\zero{g}_{01,\bar{c}}-\one{g}_{1c,\theta}),
\quad
\mone{\Gamma}{}\!^{a}{}_{0b}
=\frac{1}{2}\zero{g}{}\!^{ac}\one{g}_{bc,\theta}
-\frac{1}{2}\one{g}{}\!^{ac}\zero{g}_{bc,\theta},
\\
&&
\mone{\Gamma}{}\!^{a}{}_{bc}
=\frac{1}{2}\zero{g}{}\!^{ad}
(\zero{g}_{bd,\bar{c}}+\zero{g}_{cd,\bar{b}}-\zero{g}_{bc,\bar{d}})
+\frac{1}{2}\one{g}{}\!^{0a}\zero{g}_{bc,\theta}.
\eeaa
Nonzero connection coefficients at order $\eps^{0}$:
\beaa
&&
\zero{\Gamma}{}\!^{0}{}_{00}
=\zero{g}{}\!^{01}\zero{g}_{01,0},
\quad
\zero{\Gamma}{}\!^{0}{}_{01}
=\frac{1}{2}\one{g}{}\!^{0c}(\zero{g}_{01,\bar{c}}-\one{g}_{1c,\theta})
+\frac{1}{2}\zero{g}{}\!^{01}\two{g}_{11,\theta},
\\
&&
\zero{\Gamma}{}\!^{0}{}_{0a}
=\frac{1}{2}\zero{g}{}\!^{01}(\two{g}_{1a,\theta}+\zero{g}_{01,a})
-\frac12\one{g}{}\!^{0c}\one{g}_{ac,\theta}
+\frac12\two{h}{}\!^{0c}\zero{g}_{ac,\theta},
\\
&&
\zero{\Gamma}{}\!^{0}{}_{ab}
=\frac{1}{2}\zero{g}{}\!^{01}(\one{g}_{1b,\bar{a}}+\one{g}_{1a,\bar{b}}-\zero{g}_{ab,1})
\\
&& \qquad
-\frac12\one{g}{}\!^{0c}(\zero{g}_{bc,\bar{a}}+\zero{g}_{ac,\bar{b}}-\zero{g}_{ab,\bar{c}})
-\frac12\two{h}{}\!^{00} \zero{g}_{ab,\theta},
\\
&&
\zero{\Gamma}{}\!^{1}{}_{11}
=\frac12 \zero{g}{}\!^{01}(2\zero{g}_{01,1}-\two{g}_{11,\theta}),
\quad
\zero{\Gamma}{}\!^{1}{}_{1a}
=\frac12\zero{g}{}\!^{01}(\zero{g}_{01,a}-\two{g}_{1a,\theta}),
\\
&&
\zero{\Gamma}{}\!^{1}{}_{ab}
=-\frac12\zero{g}{}\!^{01}(\zero{g}_{ab,0}+\two{g}_{ab,\theta}),
\\
&&
\zero{\Gamma}{}\!^{a}{}_{01}
=\frac{1}{2}\zero{g}{}\!^{ac}(\two{g}_{1c,\theta}-\zero{g}_{01,c})
-\frac{1}{2}\one{g}{}\!^{ac}(\one{g}_{1c,\theta}-\zero{g}_{01,\bar{c}}),
\\
&&
\zero{\Gamma}{}\!^{a}{}_{0b}
=\frac{1}{2}\zero{g}{}\!^{ac}(\two{g}_{bc,\theta}+\zero{g}_{bc,0})
-\frac12\one{g}{}\!^{ac}\one{g}_{bc,\theta},
+\frac12\two{h}{}\!^{ac}\zero{g}_{bc,\theta},
\\
&&
\zero{\Gamma}{}\!^{a}{}_{1b}
=\frac{1}{2}\zero{g}{}\!^{ac}(\zero{g}_{bc,1}+\one{g}_{1c,\bar{b}}-\one{g}_{1b,\bar{c}})
-\frac12\one{g}{}\!^{0a}(\zero{g}_{01,\bar{b}}-\one{g}_{1b,\theta}),
\\
&&
\zero{\Gamma}{}\!^{a}{}_{bc}
= \frac12\one{g}{}\!^{0a}\one{g}_{bc,\theta}
-\frac12\two{h}{}\!^{0a}\zero{g}_{bc,\theta}
+\frac{1}{2}\zero{g}{}\!^{ad}
(\one{g}_{bd,\bar{c}}+\one{g}_{cd,\bar{b}}-\one{g}_{bc,\bar{d}})
\\
&& \qquad
-\frac{1}{2}\one{g}{}\!^{ad}
(\zero{g}_{bd,\bar{c}}+\zero{g}_{cd,\bar{b}}-\zero{g}_{bc,\bar{d}})
+\frac{1}{2}\zero{g}{}\!^{ad}
(\zero{g}_{bd,c}+\zero{g}_{cd,b}-\zero{g}_{bc,d}).
\eeaa

Nonzero component of the Ricci curvature tensor
at order $\eps^{-4}$:
\bea
&&
\mfour{R}_{00}
=-\frac{1}{2}(\zero{g}{}\!^{ab}
\zero{g}_{ab,\theta})_{,\theta}
+\frac{1}{2}\zero{g}{}\!^{01}\zero{g}_{01,\theta}
\zero{g}{}\!^{ab} \zero{g}_{ab,\theta}
-\frac{1}{4}\zero{g}{}\!^{ac} \zero{g}_{bc,\theta}
\zero{g}{}\!^{bd} \zero{g}_{ad,\theta}.
\qquad\quad
\label{Req1}
\eea

Nonzero components of the Ricci tensor at order $\eps^{-3}$:
\bea
&&
\mthree{R}_{00}
=-\frac{1}{2}\one{g}{}\!^{a}_{a,\theta\theta}
+\frac{1}{2}\zero{g}{}\!^{01}\zero{g}_{01,\theta}\one{g}{}\!^{a}_{a,\theta}
-\frac{1}{2}\zero{g}{}\!^{bd}\zero{g}_{bc,\theta}\one{g}{}\!^{c}_{d,\theta},
\label{R00-0}
\\
&&
\mthree{R}_{0a}
=\frac12(\zero{g}_{ab}\zero{g}{}\!^{01}
\one{g}{}\!^{b}_{1,\theta})_{,\theta}
+\frac14\zero{g}{}\!^{cd}\zero{g}_{cd,\theta}
\zero{g}_{ab}\zero{g}{}\!^{01} \one{g}{}\!^{b}_{1,\theta}
+\frac12(\zero{g}{}\!^{bc}\zero{g}_{ab,\theta})_{,\bar{c}}
\nonumber
\\
&&
\qquad
-\frac12(\zero{g}{}\!^{01}\zero{g}_{01,\bar{a}}
+\zero{g}{}\!^{cd}\zero{g}_{cd,\bar{a}})_{,\theta}
+\frac14\zero{g}{}\!^{bc}\zero{g}_{ab,\theta}\zero{g}{}\!^{de}\zero{g}_{de,\bar{c}}
\nonumber
\\
&&
\qquad
+\frac14\zero{g}{}\!^{01}\zero{g}_{01,\bar{a}}
\zero{g}{}\!^{cd}\zero{g}_{cd,\theta}
-\frac{1}{4}\zero{g}{}\!^{bd}\zero{g}_{cd,\theta}
\zero{g}{}\!^{ce}\zero{g}_{be,\bar{a}}.
\label{R0a-0}
\eea

The nonzero components of the Ricci curvature at
the order $\eps^{-2}$ are $\mtwo{R}_{00}$, $\mtwo{R}_{01}$, $\mtwo{R}_{0a}$ and $\mtwo{R}_{ab}$.
The components $\mtwo{R}_{00}$ and $\mtwo{R}_{0a}$ are the analogs of $\mone{R}_{0a}$ and $\mone{R}_{0a}$
for the higher order components of the metric, and will not be listed here.
The remaining components of the Ricci tensor are listed below:
\bea
&&
\mtwo{R}_{01}
=-\left(\zero{g}{}\!^{01}\zero{g}_{01,1}
+\frac12\zero{g}{}\!^{cd}\zero{g}_{cd,1}\right)_{,\theta}
-\frac14\zero{g}{}\!^{ab}\zero{g}_{bc,\theta}\zero{g}{}\!^{bd}\zero{g}_{ad,1}
\nonumber
\\
&&
+\frac12\left((\zero{g}_{01,\bar{a}}-\one{g}_{1a,\theta})\one{g}{}\!^{0a}\right)_{,\theta}
-\frac12\left((\zero{g}_{01,\bar{a}}-\one{g}_{1a,\theta})\zero{g}{}\!^{ab}\right)_{,\bar{b}}
\nonumber
\\
&&
+\frac14(\zero{g}_{01,\bar{a}}-\one{g}_{1a,\theta})
\left(\one{g}{}\!^{0a}\zero{g}{}\!^{cd}\zero{g}_{cd,\theta}
-\zero{g}{}\!^{ab}\zero{g}{}\!^{cd}\zero{g}_{cd,\bar{b}}\right),
\label{R01-1}
\\
&&
\mtwo{R}_{ab}
=-\zero{g}{}\!^{01}\zero{g}_{ab,1\theta}
+\frac{1}{2}\zero{g}{}\!^{01}\zero{g}{}\!^{cd}(\zero{g}_{ac,\theta}\zero{g}_{bd,1}
+\zero{g}_{ac,1}\zero{g}_{bd,\theta})
\nonumber
\\
&&
-\frac{1}{4}\zero{g}{}\!^{01}\zero{g}{}\!^{cd}(\zero{g}_{cd,1}\zero{g}_{ab,\theta}
+\zero{g}_{cd,\theta}\zero{g}_{ab,1})
+\stackrel{*}{R}_{ab}
\nonumber
\\
&&
+\frac{1}{2}(\one{g}{}\!^{0c}\zero{g}_{ab,\theta})_{,\bar{c}}
+\frac{1}{2}\one{g}{}\!^{0c}\zero{g}_{ab,\theta}
(\zero{g}{}\!^{01}\zero{g}_{01,\bar{c}}
+\frac12\zero{g}{}\!^{de}\zero{g}_{de,\bar{c}})
\nonumber
\\
&&
+\frac{1}{2}(\one{g}{}\!^{0c}\zero{g}_{ab,\bar{c}})_{,\theta}
+\frac{1}{2}\one{g}{}\!^{0c}\zero{g}_{ab,\bar{c}}
(\zero{g}{}\!^{01}\zero{g}_{01,\theta}
+\frac12\zero{g}{}\!^{de}\zero{g}_{de,\theta})
\nonumber
\\
&&
+\frac{1}{2}(\zero{g}{}\!^{01}\zero{g}_{ac}\one{g}{}\!^{c}_{1,\bar{b}})_{,\theta}
+\frac{1}{2}\zero{g}{}\!^{01}\zero{g}_{ac}\one{g}{}\!^{c}_{1,\bar{b}}
(\zero{g}{}\!^{01}\zero{g}_{01,\theta}
+\frac12\zero{g}{}\!^{de}\zero{g}_{de,\theta})
\nonumber
\\
&&
+\frac{1}{2}(\zero{g}{}\!^{01}\zero{g}_{bc}\one{g}{}\!^{c}_{1,\bar{a}})_{,\theta}
+\frac{1}{2}\zero{g}{}\!^{01}\zero{g}_{bc}\one{g}{}\!^{c}_{1,\bar{a}}
(\zero{g}{}\!^{01}\zero{g}_{01,\theta}
+\frac12\zero{g}{}\!^{de}\zero{g}_{de,\theta})
\nonumber
\\
&&
-\frac{1}{2}\zero{g}{}\!^{01}\zero{g}{}\!^{cd}
\zero{g}_{ac,\theta}\one{g}_{1b,\bar{d}}
+\frac{1}{2}\zero{g}{}\!^{cd}\one{g}{}\!^{0e}\zero{g}_{ac,\theta}
(\zero{g}_{be,\bar{d}}-\zero{g}_{bd,\bar{e}})
\nonumber
\\
&&
-\frac{1}{2}\zero{g}{}\!^{01}\zero{g}{}\!^{cd}
\zero{g}_{bc,\theta}\one{g}_{1a,\bar{d}}
+\frac{1}{2}\zero{g}{}\!^{cd}\one{g}{}\!^{0e}\zero{g}_{bc,\theta}
(\zero{g}_{ae,\bar{d}}-\zero{g}_{ad,\bar{e}})
\nonumber
\\
&&
-\frac{1}{2}(\one{g}{}\!^{0}_{c}\one{g}{}\!^{0c}\zero{g}_{ab,\theta})_{,\theta}
-\one{g}{}\!^{0}_{c}\one{g}{}\!^{0c}\zero{g}_{ab,\theta}
(\zero{g}{}\!^{01}\zero{g}_{01,\theta}
+\frac12\zero{g}{}\!^{de}\zero{g}_{de,\theta})
\nonumber
\\
&&
+\frac{1}{2}\one{g}{}\!^{0}_{e}\one{g}{}\!^{0e}
\zero{g}{}\!^{cd}\zero{g}_{ac,\theta}\zero{g}_{bd,\theta}
-\frac12\zero{g}{}\!^{01}\zero{g}_{ac}\one{g}{}\!^{c}_{1,\theta}
\zero{g}{}\!^{01}\zero{g}_{bd}\one{g}{}\!^{d}_{1,\theta}.
\label{Rab-1}
\eea
Here, $\stackrel{*}{R}_{ab}$ are the components of the Ricci tensor of
the zero order metric $\zero{g}_{\alpha\beta}$ regarded as a function of the
variables $\eta^a$ only,
\beaa
&&
\stackrel{*}{R}_{ab}
=\frac12\left(\zero{g}{}\!^{cd}(\zero{g}_{bd,\bar{a}}+\zero{g}_{ad,\bar{b}}
-\zero{g}_{ab,\bar{d}})\right)_{,\bar{c}}
\nonumber
\\
&&
+\frac12\zero{g}{}\!^{cd}(\zero{g}_{bd,\bar{a}}+\zero{g}_{ad,\bar{b}}-\zero{g}_{ab,\bar{d}})
(\zero{g}{}\!^{01}\zero{g}_{01,\bar{c}}+\frac12\zero{g}{}\!^{ef}\zero{g}_{ef,\bar{c}})
\nonumber
\\
&&
-(\zero{g}{}\!^{01}\zero{g}_{01,\bar{a}}+\frac12\zero{g}{}\!^{cd}\zero{g}_{cd,\bar{a}})_{,\bar{b}}
-\frac{1}{2}\zero{g}{}\!^{01}\zero{g}_{01,\bar{a}}\zero{g}{}\!^{01}\zero{g}_{01,\bar{b}}
\nonumber
\\
&&
-\frac{1}{4}\zero{g}{}\!^{ce}(\zero{g}_{de,\bar{a}}+\zero{g}_{ae,\bar{d}}
-\zero{g}_{ad,\bar{e}})\zero{g}{}\!^{df}(\zero{g}_{bf,\bar{c}}+\zero{g}_{cf,\bar{b}}
-\zero{g}_{cb,\bar{f}}).
\eeaa

\section{Variational principle}
\label{appD}
\setcounter{equation}{0}

The variational principle for the vacuum Einstein field equations is
\bea
&&\delta S = 0,\qquad S = \int L\, d^4x,\nonumber\\
&&L = R\sqrt{-\det g} ,\label{lagrangian}
\eea
where $R$ is the scalar curvature,
\[
R=g^{\alpha\beta}R_{\alpha\beta}.
\]

Using \eq{exp1}, \eq{exp2}, and \eq{Rexp} to expand the
scalar curvature, we obtain
\bea
&&R=\frac{1}{\eps^4}\mfour{R}
+\frac{1}{\eps^3}\mthree{R}
+\frac{1}{\eps^2}\mtwo{R}+O(\eps^{-1}),
\nonumber\\
&&\mfour{R}
=
\zero{g}{}\!^{\alpha\beta}\mfour{R}_{\alpha\beta},
\label{R}
\\
&&
\mthree{R}
=
\zero{g}{}\!^{\alpha\beta}\mthree{R}_{\alpha\beta}
-\one{g}{}\!^{\alpha\beta}\mfour{R}_{\alpha\beta},
\nonumber
\\
&&
\mtwo{R}
=
\zero{g}{}\!^{\alpha\beta}\mtwo{R}_{\alpha\beta}
-\one{g}{}\!^{\alpha\beta}\mthree{R}_{\alpha\beta}
+\two{h}{}\!^{\alpha\beta}\mfour{R}_{\alpha\beta}.
\nonumber
\eea
For a metric of the form \eq{rosen}, \eq{lineernst}, with $\one{g}_{11}=T$,
we find that
\bea
&&\mfour{R} = \mthree{R} = 0,\nonumber\\
&&\mtwo{R}
=2\zero{g}{}\!^{01}\mtwo{R}_{01}
+\zero{g}{}\!^{ab}\mtwo{R}_{ab} -2\one{g}\!^{0a}\mthree{R}_{0a}
+\two{h}\!^{00}\mfour{R}_{00}.
\label{scalarexp}
\eea

We use \eq{scalarexp} in \eq{lagrangian} and expand the result
with respect to $\eps$.
This gives
\[
L = \frac{1}{\eps^2} \Lmtwo + O(\eps^{-1}),
\]
with
\bea
&&
\Lmtwo
=
\left\{
2\mtwo{R}_{01}
+\zero{g}_{01}\zero{g}{}\!^{ab}\mtwo{R}_{ab}
-2\one{g}\!^{a}_1\mthree{R}_{0a} +\two{h}\!^{0}_{1}\mfour{R}_{00}
\right\}\sqrt{\det \zero{g}_{ab}}.
\qquad\quad
\label{I-2}
\eea
We make the change of variables in the integration
\[
d^4x = du\, dv\, dy\, dz = \eps^3 d\theta\, dv \, d\eta \, dz,
\]
and omit the integration with respect to the parametric variable
$z$. The leading order asymptotic variational
principle then becomes
\be
\delta \Sone = 0,\qquad \Sone = \int \Lmtwo\,d\theta\, dv \,d\eta.
\label{varprin}
\ee
Variations of $\Sone$ with respect
to the second order metric component $\two{g}_{11}$
give the constraint \eq{thetaconstraint}.
Variations with respect to the first order metric components $\one{g}_{0a}$
give the equation \eq{eq0a}.
Variations with respect to $\zero{g}_{01}$ and $\zero{g}_{ab}$
give the evolution equations \eq{apc1}.

\bigskip\bigskip\noindent
{\bf Acknowledgments.}  The work of GA was partially supported by
the CNR short-term mobility program.
The work of JKH was partially supported by the
NSF under grant number DMS--0309648.


\begin{thebibliography}{99}

\bibitem{abh} Al\`i, G., Bini, D., and Hunter, J. K., Space-times that are stationary
with respect to near light-like observers, in preparation

\bibitem{AH} Al\`i, G., and Hunter, J. K., Large Amplitude Gravitational Waves,
\emph{J. Math. Phys.}, \textbf{40} (1999), 3035--3052.

\bibitem{Br}
Brinkmann, M. W.,
On Riemann spaces conformal to Euclidean space,
{\it Proc. Nat. Acad. Sci. U.S.A.}, {\bf 9} (1923) 1--3.


\bibitem{CB}
Choquet-Bruhat, Y.,
Construction de solutions radiatives
approch\'ees des \'equations d'Einstein,
{\it Commun. Math. Phys.}, {\bf 12} (1969) 16--35.

\bibitem{CB1} Choquet-Bruhat, Y., Ondes asymptotique et approch\'{e}es
pour syst\`{e}mes nonlineaires d'\'{e}quations aux d\'{e}riv\'{e}es partielles nonlin\'{e}aires,
{\em J.~Math.~Pure et Appl.}, \textbf{48} (1969), 117-158.

\bibitem{CB2} Choquet-Bruhat, Y., The null condition and asymptotic expansions
for the Einstein equations, \emph{Ann. Phys.}, \textbf{9} (2000), 258--266.

\bibitem{CB3} Choquet-Bruhat, Y.,
Asymptotic solutions of non linear wave equations and polarized null conditions,
in \emph{Actes des Journ\'ees Math\'ematiques \`a la M\'emoire de Jean Leray}, 125--141,
S\'emin. Congr., \textbf{9}, Soc. Math. France, Paris, 2004.

\bibitem{CK} Christodoulou, D., and Klainerman, S.,
\emph{The global nonlinear stability of the Minkowski space},
Princeton Mathematical Series, \textbf{41}, Princeton University Press, Princeton, NJ, 1993.

\bibitem{GHZ} Glassey, R., Hunter, J. K., and Zheng, Y.,
 and oscillations in a nonlinear variational wave equation,
in \emph{Singularities and oscillations}, 37--60, IMA Vol. Math. Appl., \textbf{91}, Springer, New York, 1997.

\bibitem{Gr}
Griffiths, J. B., 1991,
{\it Colliding Plane Waves in General Relativity},
(Oxford University Press: Oxford).

\bibitem{Hu} Hunter, J. K., Transverse diffraction of nonlinear waves and singular rays,
{\em SIAM J.~Appl.~Math}, \textbf{48} (1988), 1-37.

\bibitem{Hrev} Hunter, J. K., Asymptotic equations for nonlinear hyperbolic waves,
in \emph{Surveys in Applied Mathematics}, \textbf{2}, 167--276, Plenum, New York, 1995.

\bibitem{HK} Hunter, J. K., and Keller, J. B., Weakly nonlinear high frequency waves,
{\em Comm.~Pure Appl.~Math.}, \textbf{36} (1983), 547-569.

\bibitem{HS} Hunter, J. K., and Saxton, R.,
Dynamics of director fields, {\em SIAM J.~Appl.~Math.}, \textbf{51} (1991), 1498-1521.

\bibitem{Is}
Isaacson, R. A.,
Gravitational radiation in the limit of high frequency I:
The linear approximation and geometrical optics
{\it Phys. Rev.}, {\bf 166} (1968), 1263--1271.

\bibitem{KP}
Khan, K. A., and Penrose, R.,
Scattering of two impulsive gravitational plane waves
{\it Nature}, {\bf 229} (1971), 185--6.

\bibitem{klainerman} Klainerman, S.,
The null condition and global existence to nonlinear wave equations,
in \emph{Nonlinear systems of partial differential equations in applied mathematics}, 293--326,
Lectures in Appl. Math., \textbf{23}, Amer. Math. Soc., Providence, RI, 1986.

\bibitem{nicolo} Klainerman, S., and Nicolo, F.,
\emph{The Evolution Problem in General Relativity},
Birkh\"auser, 2003.

\bibitem{LL}
Landau, L. D., and Lifschitz, E. M.,
{\it The Classical Theory of Fields},
Pergamon, Oxford, 1975.

\bibitem{lax} Lax, P., \emph{Hyperbolic Systems of Conservation
Laws and the Mathematical Theory of Shock Waves},
SIAM, Philadelphia, 1973.

\bibitem{lindblad} Lindblad, H., and Rodnianski, I., The weak null condition
for Einstein's equations, \emph{C. R. Acad. Sci. Paris, Ser. I }, \textbf{336}
(2003), 901--906.

\bibitem{lindblad1} Lindblad, H., and Rodnianski, I.,
Global existence for the Einstein vacuum equations in wave coordinates,
\emph{Comm. Math. Phys.}, \textbf{256} (2005), 43--110.

\bibitem{serre1} Serre, D., Oscillations non-lin\'eaires de haute fréquence; $\dim=1$,
in \emph{Nonlinear Partial Differential Equations and their Applications. Coll\`ege de France Seminar, Vol. XII},
190--210, Pitman Res. Notes Math. Ser., \textbf{302}, Longman Sci. Tech., Harlow, 1994.

\bibitem{serre2} Serre, D., Oscillations non-lin\'eaires hyperboliques de grande amplitude; $\dim\geq 2$,
in \emph{Nonlinear Variational Problems and Partial Differential Equations}, 245--294,
Pitman Res. Notes Math. Ser., \textbf{320}, Longman Sci. Tech., Harlow, 1995.

\bibitem{Sz1}
Szekeres, P.,
Colliding gravitational waves,
{\it Nature}, {\bf 228} (1970), 1183--4.

\bibitem{Sz2}
Szekeres, P., 1972,
Colliding plane gravitational waves,
{\it J. Math. Phys.}, {\bf 13} 286--94.

\bibitem{parabolic} Tappert, F. D., The parabolic approximation method,
in \emph{Wave Propagation and Underwater Acoustics}, 224--287,
Lecture Notes in Phys., \textbf{70}, Springer, Berlin, 1977.

\bibitem{taub} Taub, A. H., High-frequency gravitational waves, two-timing, and averaged Lagrangians,
in \emph{General Relativity and Gravitation}, \textbf{1}, 539--555, Plenum, New York, 1980.

\bibitem{Wh} G. B. Whitham, \emph{Linear and Nonlinear Waves}, Wiley, New York, 1974.

\end{thebibliography}
\end{document}